\documentclass[12pt,oneside,reqno]{amsart}
\usepackage{graphicx}
\usepackage{mathrsfs}
\usepackage{stmaryrd}
\usepackage{amsfonts}
\usepackage{cite}
\usepackage{enumerate,amsmath,amssymb,amsthm}
\usepackage{booktabs} 
\usepackage{diagbox} 

\newcommand{\dif}{\mathrm{d}}

\newcommand{\be}{\begin{eqnarray}}
	\newcommand{\ee}{\end{eqnarray}}
\newcommand{\ce}{\begin{eqnarray*}}
	\newcommand{\de}{\end{eqnarray*}}
\newtheorem{theorem}{Theorem}[section]
\newtheorem{lemma}[theorem]{Lemma}
\newtheorem{remark}[theorem]{Remark}
\newtheorem{definition}[theorem]{Definition}
\newtheorem{proposition}[theorem]{Proposition}
\newtheorem{Examples}[theorem]{Examples}
\newtheorem{corollary}[theorem]{Corollary}
\newtheorem{condition}[theorem]{Condition}
\def\e{\varepsilon}
\def\t{\theta}
\def\a{\alpha}

\def\d{\delta}
\def\p{\partial}

\def\s{\sigma}

\def\[{{\Big[}}
\def\]{{\Big]}}
\def\<{{\langle}}
\def\>{{\rangle}}
\def\({{\Big(}}
\def\){{\Big)}}

\def\no{\nonumber}
\def\bt{\begin{theorem}}
	\def\et{\end{theorem}}
\def\bl{\begin{lemma}}
	\def\el{\end{lemma}}
\def\br{\begin{remark}}
	\def\er{\end{remark}}
\def\bx{\begin{Examples}}
	\def\ex{\end{Examples}}
\def\bd{\begin{definition}}
	\def\ed{\end{definition}}
\def\bp{\begin{proposition}}
	\def\ep{\end{proposition}}
\def\bc{\begin{corollary}}
	\def\ec{\end{corollary}}
\def\bco{\begin{condition}}
	\def\eco{\end{condition}}
\def\cA{{\mathcal A}}

\def\cG{{\mathcal G}}

\def\cI{{\mathcal I}}
\def\cJ{{\mathcal J}}
\def\cK{{\mathcal K}}
\def\cL{{\mathcal L}}

\def\cZ{{\mathcal Z}}

\def\mE{{\mathbb E}}

\def\mH{{\mathbb H}}

\def\mN{{\mathbb N}}

\def\mP{{\mathbb P}}

\def\mR{{\mathbb R}}

\def\mV{{\mathbb V}}

\def\mX{{\mathbb X}}

\def\sA{{\mathscr A}}
\def\sB{{\mathscr B}}

\def\sF{{\mathscr F}}

\def\geq{\geqslant}
\def\leq{\leqslant}

\begin{document}
	
\allowdisplaybreaks
\title{Uniform large deviation principles and averaging principles for stochastic Burgers type equations with reflection}
	
\author{Huijie Qiao}

	
	
\thanks{This work was supported by NSF of China (No.12071071) and the Jiangsu Provincial Scientific Research Center of Applied Mathematics (No. BK20233002).}
	
	
\date{}
	
\dedicatory{School of Mathematics,
		Southeast University\\
		Nanjing, Jiangsu 211189, China\\
		hjqiaogean@seu.edu.cn}
	
\begin{abstract}
This work concerns about stochastic Burgers type equations with reflection. First of all, by means of the equicontinuous uniform Laplace principle, we prove the Freidlin-Wentzell uniform large deviation principle for these equations uniformly on bounded sets. Then based on this result, we establish the Dembo-Zeitouni uniform large deviation principle for these equations uniformly on compact sets. Finally, an averaging principle result for these equations is obtained through the time discretization approach.
\end{abstract}
	
\maketitle \rm
	
\section{Introduction}

Given a complete probability space $(\Omega, \sF, \mP)$ and a $d$-dimensional Brownian motion $W$ defined on it. Let $(\sF_t)_{t\geq 0}$ be the augmented filtration generated by $W$. Consider the following stochastic Burgers type equation with reflection:
\be\left\{\begin{array}{ll}
\dif u(t, x)=\frac{\partial^2 u(t, x)}{\partial x^2} \dif t+\frac{\partial g(t,u(t, x))}{\partial x} \dif t+f(t,x,u(t, x))\dif t+\sum\limits_{j=1}^d \sigma_j(t,x,u(t, x)) \dif W_j(t) \\
\qquad\qquad\quad+\dif K(t, x), \quad t\in\mR_+, \quad x \in[0,1],\\
u(t, x)\geq 0,\\
u(0, \cdot)=u_0(\cdot) \geq 0, \quad u(t, 0)=u(t, 1)=0,
\end{array}
\right.
\label{sbe1}
\ee
where $u_0$ is a non-negative function on $[0,1]$ and the coefficients $g: \mR_+\times\mR\to \mR$, $f: \mR_+\times[0,1]\times\mR\to\mR$, $\s_j: \mR_+\times[0,1]\times\mR\to\mR$ for $j=1,2, \cdots, d$ are Borel measurable (See Subsection \ref{sbesolu} for the definition of solutions). Stochastic partial differential equation (SPDE for short) with reflection like Eq.(\ref{sbe1}) are usually used to model the evolution of random interfaces near a hard wall (cf. \cite{fo}). Moreover, there have been a lot of results about them, such as the well-posedness (\cite{bz, dp, dz, np, xz, zhangt}), hitting properties (\cite{dmz})  and lattice approximations (\cite{zhangt0}).

Large deviation principles (LDPs for short) offer asymptotic estimates for probabilities of rare events and have been applied in many fields, such as physics, communication engineering, biology, finance and computer science. About LDPs  for SPDEs, there have been many results (See \cite{CW, cr2, dm} and references therein). Let us mention LDPs for SPDEs with reflection. Xu and Zhang \cite{xz} firstly proved the LDPs for white noise driven SPDEs with reflection by weak convergence approach. Later, Wang, Zhai and Zhang \cite{wzz}, and Brz\'ezniak, Li and Zhang \cite{blz} generalized the result in \cite{xz} to stochastic Burgers type equations with reflection and stochastic evolution equations with reflection, respectively. Besides, Matoussi, Sabbagh and Zhang \cite{msz} established the LDP for obstacle problems of quasilinear SPDEs.

Uniform large deviation principles (ULDPs for short) describe asymptotic estimates for probabilities of rare events uniformly with respect to another parameter. Moreover, they are crucial for determining the exit time and exit place of solutions. Nowadays, there have been many ULDP results for SPDEs (cf. \cite{bfz, bcf, cr2, cp, cm, eg, km, zq, ms1, sbd, ss, wang0, wang1, wang2}). For multivalued stochastic differential equations (SDEs for short), which include SDEs with reflecting on closed convex sets, Ren and Wu \cite{rw} proved a ULDP by a viscosity solution approach. Very recently, in \cite{q2} we obtained the Freidlin-Wentzell ULDP and the Dembo-Zeitouni ULDP for multivalued SDEs with jumps through the weak convergence method. However, for SPDEs with reflection, there are few related results. Therefore, the first purpose of this paper is to establish the Freidlin-Wentzell ULDP and the Dembo-Zeitouni ULDP for Eq.(\ref{sbe1}). 

Concretely speaking, consider the following stochastic Burgers type equation with reflection: 
\be\left\{\begin{array}{ll}
\dif u^\e(t, x)=\frac{\partial^2 u^\e(t, x)}{\partial x^2} \dif t+\frac{\partial g(t,u^\e(t, x))}{\partial x} \dif t+f(t,x,u^\e(t, x))\dif t \\
\qquad\qquad\quad+\sqrt\e\sum\limits_{j=1}^d \sigma_j(t,x,u^\e(t, x)) \dif W_j(t)+\dif K^\e(t, x), \quad t\in\mR_+, \quad x \in[0,1],\\
u^\e(t, x)\geq 0,\\
u^\e(0, \cdot)=u_0(\cdot) \geq 0, \quad u^\e(t, 0)=u^\e(t, 1)=0.
\end{array}
\right.
\label{sbe2}
\ee
In order to obtain the Freidlin-Wentzell ULDP for Eq.(\ref{sbe2}) uniformly on bounded sets of initial values, on the one hand, since Eq.(\ref{sbe2}) is with multiplicative noises, the technique of proving this ULDP through the uniform contraction principle is not feasible (cf. \cite{bfz, bcf, cp, cm, km, zq, wang0, wang1}). On the other hand, although the weak convergence approach is convenient (cf. \cite{bdv1, Mar, q2, jW}), it relies on the convergence of initial values, and is not straightly generalized to prove this ULDP on non-compact sets of initial values. Therefore, we show the Freidlin-Wentzell ULDP for Eq.(\ref{sbe2}) uniformly on bounded sets in terms of the equicontinuous uniform Laplace principle. Note that in \cite{ss} Salins and Setayeshgar used regularity estimates for the Green function and the equicontinuous uniform Laplace principle to obtain the Freidlin-Wentzell ULDP for Eq.(\ref{sbe2}) without reflection  in $C([0,T], L^2([0,1],\mR))$ and $C([0,T]\times[0,1])$, respectively. Here because of the reflection, there are no regularity estimates for the Green function. So, we use some new techniques to reach the goal. Moreover, our result (Theorem \ref{fwuldpsbe}) can cover \cite[Theorem 6]{ss} in one dimensional case and Theorem 3.4 in\cite{wzz}.

Besides, since the Freidlin-Wentzell ULDP and the Dembo-Zeitouni ULDP are not equivalent, we prove the Dembo-Zeitouni ULDP for Eq.(\ref{sbe2}) with the help of the compactness of level sets and the continuity of level sets in the Hausdorff metric.

Next, we observe the following stochastic Burgers type equation with reflection: 
\be\left\{\begin{array}{ll}
\dif \bar u^\e(t, x)=\frac{\partial^2 \bar u^\e(t, x)}{\partial x^2} \dif t+\frac{\partial g(t,\bar u^\e(t, x))}{\partial x} \dif t+f(\frac{t}{\e},x,\bar u^\e(t, x))\dif t \\
\qquad\qquad\quad+\sum\limits_{j=1}^d \sigma_j(\frac{t}{\e},x,\bar u^\e(t, x)) \dif W_j(t)+\dif \bar K^\e(t, x), \quad t\in\mR_+, \quad x \in[0,1],\\
\bar u^\e(t, x)\geq 0,\\
\bar u^\e(0, \cdot)=u_0(\cdot) \geq 0, \quad \bar u^\e(t, 0)=\bar u^\e(t, 1)=0.
\end{array}
\right.
\label{sbe3}
\ee
As $\e\rightarrow0$, the convergence problems of the solution $\bar u^\e$ for Eq.(\ref{sbe3}) are just averaging principles. And the averaging principle for multiscales SDEs was initiated by Khasminskii in \cite{Kh}. From then on, there are more and more averaging principle results (cf. \cite{sC, cw, dsxz, gs, my, q1, twy, yxj}). We mention some works related with ours. In \cite{yxj}, Yue, Xu and Jiao proved an averaging principle for Eq.(\ref{sbe3}) without reflection. It is natural to ask whether the averaging principle for Eq.(\ref{sbe3}) holds. Then our second purpose of this paper is to give an affirmative answer. That is, we establish that the solution $\bar u^\e$ of Eq.(\ref{sbe3}) converges to the solution $\bar u$ of Eq.(\ref{sbe4}) in probability in $C([0,T],\mH)\cap L^2([0,T],\mV)$ as $\e\rightarrow 0$. Since there are no regularity estimates for the Green function, we can not follow the line in \cite{yxj}. We utilize the penalized equation and some subtle estimations to prove this result. Note that in \cite{twy}, Tian, Wu and Yin obtained an averaging principle for stochastic evolution equations with reflection. Since the solution $\bar u^\e$ of Eq.(\ref{sbe3}) is unbounded, our result (Theorem \ref{averprintheo}) can not be covered by \cite[Theorem 3.1]{twy}.

The rest of this paper is organized as follows. In Section \ref{pre}, we introduce some notations and ULDPs. In Section \ref{main} the main results are stated. We arrange the proofs of three main results in Section \ref{fwuldpproo}, Section \ref{dzuldpproo} and Section \ref{averprintheoproo}, respectively.

The following convention will be used throughout the paper: $C$ with or without indices will denote different positive constants whose values may change from one place to another.

\section{Preliminaries}\label{pre}

In this section, we introduce some notations and ULDPs.

\subsection{Notation}

In this subsection, we introduce some notations used in the sequel.

Let $C_c([0,1])$ be the set of all continuous functions on $[0,1]$ with compact supports. Let $C^2_c([0,1])$ and  $C^\infty_c([0,1])$  be the subsets of $C_c([0,1])$ where all functions have $2$ and infinite order derivatives, respectively. 

Let $\mH=L^2([0,1], \mR)$ be the usual $L^2$-space with norm $|\cdot|_\mH$ and inner product $\<\cdot, \cdot\>$. Denote by $\mV$ the Sobolev space of order one, i.e., $\mV$ is the completion of $C_c^{\infty}([0,1])$ under the norm $\|u\|_{\mV}^2=\int_0^1\left(\frac{\partial u}{\partial x}\right)^2 d x$. The corresponding inner product will be denoted by $\<\<\cdot, \cdot\>\>$. Remark that $\mV=\left\{u \in \mH^{1}([0,1]): u(0)=u(1)=0\right\}$, where $\mH^{1}([0,1])$ denotes the usual Sobolev space of absolutely continuous functions defined on $[0,1]$ whose derivative belongs to $\mH$. $\mV^*$ denotes the dual space of $\mV$. 

\subsection{Stochastic Burgers type equations with reflection}\label{sbesolu}

In this subsection, we introduce the solutions to stochastic Burgers type equations with reflection.

We recall Eq.(\ref{sbe1}), i.e.
\ce\left\{\begin{array}{ll}
\dif u(t, x)=\frac{\partial^2 u(t, x)}{\partial x^2} \dif t+\frac{\partial g(t,u(t, x))}{\partial x} \dif t+f(t,x,u(t, x))\dif t+\sum\limits_{j=1}^d \sigma_j(t,x,u(t, x)) \dif W_j(t) \\
\qquad\qquad\quad+\dif K(t, x), \quad t\in\mR_+, \quad x \in[0,1],\\
u(t, x)\geq 0,\\
u(0, \cdot)=u_0(\cdot) \geq 0, \quad u(t, 0)=u(t, 1)=0.
\end{array}
\right.
\de
The following is the definition of a solution to Eq.(\ref{sbe1}).

\bd\label{soludefi}
A pair $(u, K)$ is said to be a solution of Eq.(\ref{sbe1}) if

$(i)$ $u(t, \cdot)$ is $\mV$-valued $\sF_t$-measurable for any $t\geq 0$ and $u(t, x) \geq 0$ a.e. for any $(t,x)\in\mR_+\times[0,1]$;

$(ii)$ $K$ is a random measure on $\mR_+\times[0,1]$ such that

$(a)$ $\mE\left[({\rm Var}(K)([0,T]\times[0,1]))^2\right]<+\infty, \forall T \geq 0$, where ${\rm Var}(K)([0,T]\times[0,1])$ denotes the total variation of $K$ on $[0, T]\times[0,1]$ defined by 
\ce
{\rm Var}(K)([0,T]\times[0,1]):=\sup\limits_{\pi}\sum_{i=1}^n |K(E_i)|,
\de
and the supremum is taken over all partitions $\pi$ of the domain $[0,T]\times[0,1]$;

$(b)$ $K$ is adapted in the sense that for any bounded measurable mapping $\phi$:
\ce
\int_0^t\int_0^1\phi(s,x)K(\dif s, \dif x) ~\mbox{is}~\sF_t\mbox{-measurable}.
\de

$(iii)$ $(u, K)$ solves the parabolic stochastic partial differential equation with reflection in the following sense: $\forall t \in \mathbb{R}_{+}, \psi \in C_c^2([0,1])$ with $\psi(0)=\psi(1)=0$,
\ce
\<u(t), \psi\>&=&\<u(0), \psi\>+\int_0^t\<u(s), \psi^{\prime \prime}\>\dif s-\int_0^t\<g(s,u(s)), \psi^{\prime}\>\dif s+\int_0^t\<f(s,\cdot,u(s)),\psi\>\dif s\\
&&+\sum_{j=1}^d \int_0^t\<\sigma_j(s,\cdot,u(s)), \psi\>\dif W_j(s)+\int_0^t\int_0^1\psi(x)K(\dif s,\dif x), \quad a.s.
\de
where $u(t):=u(t, \cdot)$;

$(iv)$ for any $T>0$, $\int_0^T\int_0^1u(t,x)K(\dif t,\dif x)=0$.
\ed

If $(u, K)$ is a solution of Eq.(\ref{sbe1}) in the sense of Definition \ref{soludefi}, for any $T>0$ and $t\in[0,T]$,
\ce
u(t)&=&u_0+\int_0^t\frac{\partial^2 u(s)}{\partial x^2} \dif s+\int_0^t\frac{\partial g(s,u(s))}{\partial x} \dif s+\int_0^t f(s,\cdot,u(s))\dif s\\
&&+\sum\limits_{j=1}^d\int_0^t \sigma_j(s,\cdot,u(s)) \dif W_j(s)+\int_0^tK(\dif s,\dif x)
\de
in $\mV^*$, $\mP$-a.s..

\subsection{Uniform large deviation principles}

In this subsection, we introduce the ULDPs (cf. \cite{ms}).

Let $(\mX, \rho_\mX)$ be a Polish space and let $\mX_0$ be a set. Let $\{u^{\e}_{u_0}, \e>0, u_0\in \mX_0\}$ be a family of $\mX$-valued random variables defined on $(\Omega, \mathscr{F}, \{\mathscr{F}_t\}_{t\geq 0}, \mP)$.

\bd\label{rfde} 
$(i)$ For $u_0\in\mX_0$, a function $\Lambda_{u_0}: \mX\rightarrow[0,+\infty]$ is called a rate function on $\mX$, if for all $M\geq 0$, $\{\phi\in\mX: \Lambda_{u_0}(\phi)\leq M\}$ is a closed subset of $\mX$.

$(ii)$ For $u_0\in\mX_0$, a function $\Lambda_{u_0}: \mX\rightarrow[0,+\infty]$ is called a good rate function on $\mX$, if for all $M\geq 0$, $\{\phi\in\mX: \Lambda_{u_0}(\phi)\leq M\}$ is a compact subset of $\mX$.
\ed

Here we recall the Freidlin-Wentzell formulation of the uniformly large deviation principle. This definition can also be found in \cite{ms}.

\bd [Freidlin-Wentzell ULDP]\label{fwuldpde}
Let $\sA$ be a collection of subsets of $\mX_0$. $\{u^{\e}_{u_0}, \e>0, u_0\in \mX_0\}$ are said to satisfy a Freidlin-Wentzell ULDP with respect to the rate function $\Lambda_{u_0}$ with the speed $\e$ uniformly over $\sA$, if

$(a)$ For any $B\in \sA, M, \d, \t>0$, there exists a $\e_0>0$ such that
\ce
\mP(\rho_\mX(u^{\e}_{u_0}, \phi)<\d)\geq \exp\left\{-\frac{\Lambda_{u_0}(\phi)+\t}{\e}\right\},
\de
for all $0<\e<\e_0, u_0\in B, \phi\in \Phi_{u_0}(M)$, where $\Phi_{u_0}(M):=\{\phi\in\mX: \Lambda_{u_0}(\phi)\leq M \}$.

$(b)$ For any $B\in \sA, M, \d, \t>0$, there exists a $\e_0>0$ such that
\ce
\mP(\rho_\mX(u^{\e}_{u_0}, \Phi_{u_0}(M'))\geq\d)\leq \exp\left\{-\frac{M'-\t}{\e}\right\},
\de
for all $0<\e<\e_0, u_0\in B, 0\leq M'\leq M$, where 
\ce
\rho_\mX(u^{\e}_{u_0}, \Phi_{u_0}(M')):=\inf\limits_{\phi\in\Phi_{u_0}(M')} \rho_\mX(u^{\e}_{u_0}, \phi).
\de
\ed

Next, in order to prove the Freidlin-Wentzell ULDP, we also introduce the equicontinuous uniform Laplace principle. A family $\cL$ of functions $\Psi: \mX \rightarrow \mathbb{R}$ is called equibounded and equicontinuous if
\ce
\sup _{\Psi \in\cL} \sup _{\phi \in \mX}|\Psi(\phi)|<\infty \text { and } \lim _{\delta \rightarrow 0} \sup _{\substack{\Psi \in \cL}} \sup _{\substack{\phi_1, \phi_2 \in \mX \\\rho_{\mX}(\phi_1, \phi_2)<\d}}|\Psi(\phi_1)-\Psi(\phi_2)|=0.
\de

\bd[Equicontinuous uniform Laplace principle]\label{ulpde} 
Let $\sA$ be a collection of subsets of $\mX_0$. $\{u^{\e}_{u_0}, \e>0, u_0\in \mX_0\}$ are said to satisfy an equicontinuous uniform Laplace principle with respect to the rate function $\Lambda_{u_0}$ with the speed $\e$ uniformly over $\sA$, if for any $B\in \sA$ and any equicontinuous and equibounded family $\cL$,
\ce
\lim\limits_{\e\rightarrow 0}\sup _{\Psi\in \cL} \sup\limits_{u_0\in B}\left|\e \log \mE \exp\left(-\frac{\Psi(u^{\e}_{u_0})}{\e}\right)+\inf\limits_{\phi\in \mX}\left\{\Psi(\phi)+\Lambda_{u_0}(\phi)\right\}\right|=0.
\de
\ed

Note that the Freidlin-Wentzell ULDP and the equicontinuous uniform Laplace principle are equivalent (\cite[Theorem 2.10]{ms}). Thus, in order to obtain the Freidlin-Wentzell ULDP, we present a sufficient condition for the equicontinuous uniform Laplace principle. Set for any $N>0$ and $T>0$, 
$$
\mathbf{D}^N=\left\{h: [0, T]\rightarrow\mR^d| h ~\mbox{is}~ \sB([0,T])/\sB(\mR^d) ~\mbox{measurable, and}~ \int_{0}^{T}|h(t)|^{2}\dif t \leq N\right\},
$$
and we equip $\mathbf{D}^{N}$ with the weak convergence topology in $L^2\left([0, T], \mR^d\right)$. So, $\mathbf{D}^{N}$ is metrizable as a compact Polish space. Let $\cA$ be  the collection of $(\sF_t)_{t\in[0,T]}$-adapted square integrable $\mathbb{R}^d$-valued processes and $\mathcal{A}^N$ be the space of $\mathbf{D}^N$-valued random controls:
$$
\cA^N=\left\{\xi\in \mathcal{A}: \xi(\cdot, \omega) \in \mathbf{D}^N, \mathbb{P}\text { a.s. } \omega\right\}.
$$

Let $\mathcal{G}^\e_{u_0}: C([0,T],\mR^d)\rightarrow \mX$ be a measurable mapping for any $0<\e<1$ and $u_0\in\mX_0$.

\bco\label{cond}
There exists a measurable mapping $\mathcal{G}^0_{u_0}: C([0,T],\mR^d)\rightarrow \mX$ such that for any $\delta>0$ and $B \in \mathscr{A}$
$$
\lim _{\varepsilon \rightarrow 0} \sup _{u_0 \in B} \sup _{\xi\in \cA^N} \mathbb{P}\left(\rho_{\mX}\left(\cG^{\varepsilon}_{u_0}\left(\sqrt{\varepsilon} W+\int_0^{.} \xi(s) \dif s\right), \cG^0_{u_0}\left(\int_0^{\cdot} \xi(s) \dif s\right)\right)>\delta\right)=0.
$$
\eco

For any $u_0\in\mX_0$ and $\phi \in \mX$, set 
$$
\mathbf{D}_\phi:=\left\{h\in L^2\left([0, T], \mR^d\right): \phi=\mathcal{G}^0_{u_0}\left(\int_0^{\cdot} h(s) \dif s\right)\right\}.
$$
Then we define
\ce
\Lambda_{u_0}(\phi)=\left\{\begin{array}{ll}
\inf\limits _{h\in\mathbf{D}_\phi}\left\{\frac{1}{2}\int_0^T|h(t)|^2\dif t\right\}, \quad \mathbf{D}_\phi\neq\emptyset,\\
\infty, \qquad\qquad\qquad\qquad\qquad\mathbf{D}_\phi=\emptyset.
\end{array}
\right.
\de

The following result is due to \cite[Theorem 2.13]{ms} and \cite[Theorem 2.10]{ms}.

\bt\label{uldpfwjuth}
Set $u^{\e}_{u_0}=\mathcal{G}^\e_{u_0}\left(\sqrt{\e} W\right)$. Suppose that Condition \ref{cond} holds. Then, the family $\left\{u^{\e}_{u_0}\right\}$ satisfies a Freidlin-Wentzell ULDP on $\mX$, with the rate function $\Lambda_{u_0}$, uniformly over $\sA$.
\et

Finally, we introduce the definition of the other uniform large deviations principle which is called the Dembo-Zeitouni ULDP. The following definition can be found in \cite{ms}.

\bd[Dembo-Zeitouni ULDP]
Let $\mathscr{A}$ be a collection of subsets of $\mX_0$. $\{u^{\e}_{u_0}, \e>0, u_0\in \mX_0\}$ are said to satisfy a Dembo-Zeitouni ULDP with respect to the rate function $\Lambda_{u_0}$ with the speed $\e$ uniformly over $\mathscr{A}$, if

$(a)$ For any $B \in \mathscr{A}$ and any open $G \subset \mX$,

$$
\liminf _{\varepsilon \rightarrow 0} \inf _{u_0\in B}\left(\e\log \mathbb{P}\left(u^{\e}_{u_0} \in G\right)\right) \geq-\sup _{u_0\in B} \Lambda_{u_0}(G)
$$

(b) For any $B \in \mathscr{A}$ and any closed $F \subset \mX$,

$$
\limsup _{\varepsilon \rightarrow 0} \sup _{u_0\in B}\left(\e \log \mathbb{P}\left(u^{\e}_{u_0}\in F\right)\right) \leq-\inf _{u_0\in B} \Lambda_{u_0}(F) .
$$
\ed

Note that according to\cite{ms}, these two definitions of ULDP are not equivalent.  Thus, in order to obtain the Dembo-Zeitouni ULDP, we assume:

\bco\label{abcd}
$(a)$ $\mX_0$ is a Polish space with metric $\rho_{\mX_0}$;

$(b)$ $\mathscr{A}$ is the collection of compact subsets of $\mX_0$;

$(c)$ For every $u_0 \in \mX_0, \Lambda_{u_0}$ is a good rate function;

$(d)$ The level sets are continuous in the Hausdorff metric in the sense that as \\$\lim\limits _{n \rightarrow\infty} \rho_{\mX_0}\left(u_{0,n}, u_0\right)=0$, for any $M \geq 0$,
$$
 \lim _{n \rightarrow+\infty} \max\left\{\sup\limits_{\phi\in\Phi_{u_0}(M)}\rho_\mX\(\phi,\Phi_{u_{0,n}}(M)\), \sup\limits_{\phi\in\Phi_{u_{0,n}}(M)}\rho_\mX\(\phi,\Phi_{u_{0}}(M)\)\right\}=0.
$$
\eco

By Theorem \ref{uldpfwjuth} and \cite[Theorem 2.7]{ms}, we obtain the following result.

\bt\label{uldpdzjuth}
Set $u^{\e}_{u_0}=\mathcal{G}^\e_{u_0}\left(\sqrt{\e} W\right)$. Suppose that Condition \ref{cond} and \ref{abcd} hold. Then, the family $\left\{u^{\e}_{u_0}\right\}$ satisfies a Dembo-Zeitouni ULDP on $\mX$, with the rate function $\Lambda_{u_0}$, uniformly over $\sA$.
\et

\section{Main results}\label{main}

In this section, we formulate the main results in this paper. 

\subsection{Two uniformly large deviation principles for stochastic Burgers type equations with reflection} 

In this subsection, we state the Freidlin-Wentzell ULDP and the Dembo-Zeitouni ULDP results for stochastic Burgers type equations with reflection.

We recall Eq.(\ref{sbe2}), i.e.
\ce\left\{\begin{array}{ll}
\dif u^\e(t, x)=\frac{\partial^2 u^\e(t, x)}{\partial x^2} \dif t+\frac{\partial g(t,u^\e(t, x))}{\partial x} \dif t+f(t,x,u^\e(t, x))\dif t \\
\qquad\qquad\quad+\sqrt\e\sum\limits_{j=1}^d \sigma_j(t,x,u^\e(t, x)) \dif W_j(t)+\dif K^\e(t, x), \quad t\in\mR_+, \quad x \in[0,1],\\
u^\e(t, x)\geq 0,\\
u^\e(0, \cdot)=u_0(\cdot) \geq 0, \quad u^\e(t, 0)=u^\e(t, 1)=0.
\end{array}
\right.
\de

Assume:
\begin{enumerate}[$(\mathbf{H}_{g})$]
\item $g(s,z)$ is continuously differentiable in $z$ and for any $T>0$ there exists a constant $L_g>0$ such that for any $s\in[0,T]$ and $z\in\mR$
\ce
\left|\frac{\p g(s,z)}{\p z}\right|\leq L_g(1+|z|).
\de
\end{enumerate}
\begin{enumerate}[$(\mathbf{H}_{f})$]
\item $f(s,x,z)$ is continuous in $z$ and for any $T>0$, there exists a constant $L_f>0$ such that for any $s\in[0,T]$, $x\in[0,1]$ and $z_1, z_2, z\in \mR$
\ce
&&(z_1-z_2)(f(s,x,z_1)-f(s,x,z_2))\leq L_f|z_1-z_2|^2,\\
&&|f(s,x,z)|^2\leq L_f(1+|z|^2).
\de
\end{enumerate}
\begin{enumerate}[$(\mathbf{H}_{\s})$]
\item For any $T>0$, there exists a constant $L_\s>0$ such that for any $s\in[0,T]$, $x\in[0,1]$ and $z_1, z_2, z\in \mR$
\ce
&&\sum\limits_{j=1}^d|\s_j(s,x,z_1)-\s_j(s,x,z_2)|^2\leq L_\s|z_1-z_2|^2,\\
&&\sum\limits_{j=1}^d|\s_j(s,x,z)|^2\leq L_\s(1+|z|^2).
\de
\end{enumerate}

\br
$(\mathbf{H}_{g})$ implies that for any $z_1, z_2\in \mR$
\ce
|g(s,z_1)-g(s,z_2)|\leq L_g(1+|z_1|+|z_2|)|z_1-z_2|.
\de
\er

By the similar deduction to that for \cite[Theorem 3.1]{zhangt}, we obtain the following well-posedness result for Eq.(\ref{sbe2}).

\bt\label{wellpose}
Assume that $(\mathbf{H}_{g})$, $(\mathbf{H}_{f})$ and $(\mathbf{H}_{\s})$ hold and $u_0\in\mV$. Then Eq.(\ref{sbe2}) has a unique solution $(u^\e_{u_0},K^\e_{u_0})$ that satisfies
\ce
\mE\left[\sup\limits_{t\in[0,T]}|u^\e_{u_0}(t)|_\mH^2\right]+\mE\left[\int_0^T\|u^\e_{u_0}(t)\|_\mV^2\dif t\right]<\infty, \quad \forall T>0.
\de
\et

Next, for any $T>0$ and $h\in L^2([0,T],\mR^d)$, consider the following Burgers type equation with reflection
\be
\left\{\begin{array}{ll}
\dif u^{0,h}(t, x)=\frac{\partial^2 u^{0,h}(t, x)}{\partial x^2} \dif t+\frac{\partial g(t,u^{0,h}(t, x))}{\partial x} \dif t+f(t,x,u^{0,h}(t, x))\dif t\\
\qquad\qquad\qquad+\sum\limits_{j=1}^d \sigma_j(t,x,u^{0,h}(t, x))h_j(t) \dif t+\dif K^{0,h}(t, x), \quad t\in[0,T], \quad x \in[0,1],\\
u^{0,h}(t, x)\geq 0,\\
u^{0,h}(0, \cdot)=u_0(\cdot) \geq 0, \quad u^{0,h}(t, 0)=u^{0,h}(t, 1)=0.
\end{array}
\right.
\label{conequ0h}
\ee
Under $(\mathbf{H}_{g})$, $(\mathbf{H}_{f})$ and $(\mathbf{H}_{\s})$, by the similar deduction to that for Theorem \ref{wellpose}, we know that Eq.(\ref{conequ0h}) has a unique 
solution $(u^{0,h}_{u_0},K^{0,h}_{u_0})$ with
$$
\sup\limits_{t\in[0,T]}|u^{0,h}_{u_0}(t)|_\mH^2+\int_0^T\|u^{0,h}_{u_0}(t)\|_\mV^2\dif t<\infty.
$$
For any $\phi\in C([0,T],\mH)\cap L^2([0,T],\mV)$, define $\Lambda_{u_0}: C([0,T],\mH)\cap L^2([0,T],\mV)\rightarrow [0,\infty]$ by
\be
\Lambda_{u_0}(\phi):=\frac{1}{2}\inf\limits_{\phi=u^{0,h}_{u_0}}\int_0^T|h(t)|^2\dif t,
\label{ratefunc}
\ee
and $\Lambda_{u_0}$ is a rate function.

Now, we state the Freidlin-Wentzell ULDP for Eq.(\ref{sbe1}), which is the first main result in this subsection.

\bt\label{fwuldpsbe}
Assume that $(\mathbf{H}_{g})$, $(\mathbf{H}_{f})$ and $(\mathbf{H}_{\s})$ hold. Then the family $\{u^\e_{u_0}\}$ satisfies the Freidlin-Wentzell ULDP on $C([0,T],\mH)\cap L^2([0,T],\mV)$ with the rate function $\Lambda_{u_0}$, uniformly on bounded subsets of $\mV$.
\et

\br
We mention that our Freidlin-Wentzell ULDP for Eq.(\ref{sbe1}) holds uniformly on bounded subsets of $\mV$ rather than on compact subsets of $\mV$.
\er

The proof of Theorem \ref{fwuldpsbe} is postponed to Section \ref{fwuldpproo}.

Then we present the Dembo-Zeitouni ULDP for Eq.(\ref{sbe1}), which is the second main result in this subsection.

\bt\label{dzuldpsbe}
Assume that $(\mathbf{H}_{g})$, $(\mathbf{H}_{f})$ and $(\mathbf{H}_{\s})$ hold. Then the family $\{u^\e_{u_0}\}$ satisfies the Dembo-Zeitouni ULDP on $C([0,T],\mH)\cap L^2([0,T],\mV)$ with the good rate function $\Lambda_{u_0}$, uniformly on compact subsets of $\mV$.
\et

The proof of Theorem \ref{dzuldpsbe} is placed in Section \ref{dzuldpproo}.

\subsection{An averaging principle for stochastic Burgers type equations with reflection}

In this subsection, we present an averaging principle result for stochastic Burgers type equations with reflection.

Consider Eq.(\ref{sbe3}), i.e.
\ce\left\{\begin{array}{ll}
\dif \bar u^\e(t, x)=\frac{\partial^2 \bar u^\e(t, x)}{\partial x^2} \dif t+\frac{\partial g(t,\bar u^\e(t, x))}{\partial x} \dif t+f(\frac{t}{\e},x,\bar u^\e(t, x))\dif t \\
\qquad\qquad\quad+\sum\limits_{j=1}^d \sigma_j(\frac{t}{\e},x,\bar u^\e(t, x)) \dif W_j(t)+\dif \bar K^\e(t, x), \quad t\in\mR_+, \quad x \in[0,1],\\
\bar u^\e(t, x)\geq 0,\\
\bar u^\e(0, \cdot)=u_0(\cdot) \geq 0, \quad \bar u^\e(t, 0)=\bar u^\e(t, 1)=0.
\end{array}
\right.
\de
Under $(\mathbf{H}_{g})$, $(\mathbf{H}_{f})$ and $(\mathbf{H}_{\s})$, by Theorem \ref{wellpose}, we know that for $u_0\in\mV$, Eq.(\ref{sbe3}) has a unique 
solution $(\bar u^{\e}_{u_0},\bar K^{\e}_{u_0})$ with
$$
\mE\left[\sup\limits_{t\in[0,T]}|\bar u^{\e}_{u_0}(t)|_\mH^2+\int_0^T\|\bar u^{\e}_{u_0}(t)\|_\mV^2\dif t\right]<\infty, \quad \forall T>0.
$$
In order to study the convergence of $\bar u^{\e}_{u_0}$ as $\e\rightarrow 0$, we also assume:
\begin{enumerate}[$(\mathbf{H}^C_{g})$]
\item $g(s,z)$ is continuous in $s$.
\end{enumerate}
\begin{enumerate}[$(\mathbf{H}^L_{f,\s})$]
\item
There exist $\bar{f}: [0,1]\times\mR\rightarrow\mR$ and $\bar \s_j: [0,1]\times\mR\rightarrow\mR$ for $j=1,2, \cdots, d$ satisfying for any $\hat{T}\in\mR_+, x\in[0,1], z\in\mR$
\ce
&&\frac{1}{\hat{T}}\int_0^{\hat{T}}\left|f(s,x,z)-\bar{f}(x,z)\right|^2\dif s\leq \kappa(\hat{T})(1+|z|^2),\\
&&\frac{1}{\hat{T}}\int_0^{\hat{T}}\sum\limits_{j=1}^d\left|\s_j(s,x,z)-\bar{\s}_j(x,z)\right|^2\dif s\leq \kappa(\hat{T})(1+|z|^2),
\de
where $\kappa(\cdot)$ is a continuous and positive bounded function with $\lim\limits_{\hat{T}\rightarrow\infty}\kappa(\hat{T})=0$.
\end{enumerate}

\br\label{barflip}
$(i)$ By $(\mathbf{H}_{g})$ and $(\mathbf{H}^C_{g})$, it holds that for any $T>0$ and $s\in[0,T]$
\ce
|g(s,z)|\leq \sup\limits_{s\in[0,T]}|g(s,0)|+L_g(1+|z|)|z|.
\de
That is, $g(s,z)$ is of polynomial growth in $z$.

$(ii)$ By $(\mathbf{H}_{f})$, $(\mathbf{H}_{\s})$ and $(\mathbf{H}^L_{f,\s})$, we know that for any $x\in[0,1]$ and $z_1, z_2, z\in \mR$
\ce
&&(z_1-z_2)(\bar f(x,z_1)-\bar f(x,z_2))\leq L_f|z_1-z_2|^2,\\
&&|\bar f(x,z)|^2\leq 2L_f(1+|z|^2),\\
&&\sum\limits_{j=1}^d|\bar\s_j(x,z_1)-\bar\s_j(x,z_2)|^2\leq 3L_\s|z_1-z_2|^2,\\
&&\sum\limits_{j=1}^d|\bar\s_j(x,z)|^2\leq 2L_\s(1+|z|^2).
\de
Since the verifications of four above inequalities are similar, we only verify the first inequality. Indeed, for $z_1, z_2\in\mR$,
\ce
&&(z_1-z_2)(\bar{f}(x, z_1)-\bar{f}(x, z_2))\\
&=&(z_1-z_2)\left(\bar{f}(x, z_1)-\frac{1}{\hat{T}}\int_0^{\hat{T}}f(s,x,z_1)\dif s\right)\\
&&+(z_1-z_2)\left(\frac{1}{\hat{T}}\int_0^{\hat{T}}f(s,x,z_1)\dif s-\frac{1}{\hat{T}}\int_0^{\hat{T}}f(s,x,z_2)\dif s\right)\\
&&+(z_1-z_2)\left(\frac{1}{\hat{T}}\int_0^{\hat{T}}f(s,x,z_2)\dif s-\bar{f}(x,z_2)\right)\\
&\leq&|z_1-z_2|\frac{1}{\hat{T}}\int_0^{\hat{T}}\left|f(s,x,z_1)-\bar{f}(x,z_1)\right|\dif s+\frac{1}{\hat{T}}\int_0^{\hat{T}}(z_1-z_2)(f(s,x,z_1)-f(s,x,z_2))\dif s\\
&&+|z_1-z_2|\frac{1}{\hat{T}}\int_0^{\hat{T}}\left|f(s,x,z_2)-\bar{f}(x,z_2)\right|\dif s\\
&\leq&|z_1-z_2|\left(\frac{1}{\hat{T}}\int_0^{\hat{T}}\left|f(s,x,z_1)-\bar{f}(x,z_1)\right|^2\dif s\right)^{1/2}+L_f|z_1-z_2|^2\\
&&+|z_1-z_2|\left(\frac{1}{\hat{T}}\int_0^{\hat{T}}\left|f(s,x,z_2)-\bar{f}(x,z_2)\right|^2\dif s\right)^{1/2}\\
&\leq&|z_1-z_2|\kappa^{1/2}(\hat{T})(2+|z_1|+|z_2|)+L_f|z_1-z_2|^2.
\de
Let $\hat{T}\rightarrow\infty$, and we obtain that $(z_1-z_2)(\bar{f}(x, z_1)-\bar{f}(x, z_2))\leq L_f|z_1-z_2|^2$. 
\er

Next, we construct the following averaging stochastic Burgers type equations with reflection:
\be\left\{\begin{array}{ll}
\dif \bar u^0(t, x)=\frac{\partial^2 \bar u^0(t, x)}{\partial x^2} \dif t+\frac{\partial g(t,\bar u^0(t, x))}{\partial x} \dif t+\bar f(x,\bar u^0(t, x))\dif t \\
\qquad\qquad\quad+\sum\limits_{j=1}^d \bar\sigma_j(x,\bar u^0(t, x)) \dif W_j(t)+\dif \bar K^0(t, x), \quad t\in\mR_+, \quad x \in[0,1],\\
\bar u^0(t, x)\geq 0,\\
\bar u^0(0, \cdot)=u_0(\cdot) \geq 0, \quad \bar u^0(t, 0)=\bar u^0(t, 1)=0.
\end{array}
\right.
\label{sbe4}
\ee
Under $(\mathbf{H}_{g})$, $(\mathbf{H}_{f})$, $(\mathbf{H}_{\s})$ and $(\mathbf{H}^L_{f,\s})$, by Remark \ref{barflip} and Theorem \ref{wellpose}, we know that for $u_0\in\mV$, Eq.(\ref{sbe4}) has a unique solution $(\bar u^{0}_{u_0},\bar K^{0}_{u_0})$ with
$$
\mE\left[\sup\limits_{t\in[0,T]}|\bar u^{0}_{u_0}(t)|_\mH^2+\int_0^T\|\bar u^{0}_{u_0}(t)\|_\mV^2\dif t\right]<\infty, \quad \forall T>0.
$$
Then the following theorem describes the relationship between $\bar u^{\e}_{u_0}$ and $\bar u^{0}_{u_0}$.

\bt\label{averprintheo}
Assume that $(\mathbf{H}_{g})$, $(\mathbf{H}_{f})$, $(\mathbf{H}_{\s})$, $(\mathbf{H}^C_{g})$ and $(\mathbf{H}^L_{f,\s})$ hold and $u_0\in\mV$. Then for any $T>0$, $\bar u^{\e}_{u_0}$ converges to $\bar u^{0}_{u_0}$ in probability in $C([0,T],\mH)\cap L^2([0,T],\mV)$ as $\e\rightarrow 0$.
\et

The proof of Theorem \ref{averprintheo} is presented in Section \ref{averprintheoproo}.

\section{Proof of Theorem \ref{fwuldpsbe}}\label{fwuldpproo}

In this section, we prove Theorem \ref{fwuldpsbe}. We first prepare some estimates, verify Condition \ref{cond} and then show Theorem \ref{fwuldpsbe} by Theorem \ref{uldpfwjuth}.

First of all, the unique solution to Eq.(\ref{sbe2}) is a strong solution in the sense of probability theory. Therefore, for any $\e>0$ and $u_0\in\mH$, there exists a measurable mapping $\mathcal{G}^\e_{u_0}: C([0,T],\mR^d)\rightarrow C([0,T],\mH)\cap L^2([0,T],\mV)$ such that $u^\e_{u_0}=\mathcal{G}^\e_{u_0}\left(\sqrt{\e} W\right)$.

For any $\xi\in\mathcal{A}^N$, consider the following stochastic Burgers type equation with reflection:
\be\left\{\begin{array}{ll}
\dif u^{\e,\xi}(t, x)=\frac{\partial^2 u^{\e,\xi}(t, x)}{\partial x^2} \dif t+\frac{\partial g(t,u^{\e,\xi}(t, x))}{\partial x} \dif t+f(t,x,u^{\e,\xi}(t, x))\dif t\\
\qquad\qquad\qquad+\sqrt\e\sum\limits_{j=1}^d \sigma_j(t,x,u^{\e,\xi}(t, x)) \dif W_j(t) +\sum\limits_{j=1}^d \sigma_j(t,x,u^{\e,\xi}(t, x)) \xi(s)\dif s\\
\qquad\qquad\qquad+\dif K^{\e,\xi}(t, x), \quad t\in[0,T],\quad x \in[0,1],\\
u^{\e,\xi}(t, x)\geq 0,\\
u^{\e,\xi}(0, \cdot)=u_0(\cdot) \geq 0, \quad u^{\e,\xi}(t, 0)=u^{\e,\xi}(t, 1)=0.
\end{array}
\right.
\label{conequxih}
\ee
By the Girsanov theorem, we know that the above equation has a unique solution $(u^{\e,\xi}_{u_0},\\ K^{\e,\xi}_{u_0})$ and $u^{\e,\xi}_{u_0}=\mathcal{G}^\e_{u_0}\left(\sqrt{\e} W+\int_0^{.} \xi(s) \dif s\right)$.

Next, we define the other measurable mapping $\mathcal{G}^0_{u_0}: C([0,T],\mR^d)\rightarrow C([0,T],\mH)\cap L^2([0,T],\mV)$ by
\ce
u^{0,\xi}_{u_0}=\cG^0_{u_0}(\int_0^{.} \xi(s) \dif s),
\de
where $u^{0,\xi}_{u_0}$ solves Eq.(\ref{conequ0h}) by replacing $h$ by $\xi$. We estimate $u^{\e,\xi}_{u_0}$ and $u^{0,\xi}_{u_0}$.

\bl\label{uesefomoes}
Under $(\mathbf{H}_{g})$, $(\mathbf{H}_{f})$ and $(\mathbf{H}_{\s})$, there exists a constant $C>0$ such that
\be
&&\sup\limits_{\e\in(0,1)}\left\{\mE\sup\limits_{t\in[0,T]}|u^{\e,\xi}_{u_0}(t)|^2_\mH+\mE\int_0^T\|u^{\e,\xi}_{u_0}(t)\|^2_\mV\dif t\right\}\leq C(1+|u_0|^2_\mH),\label{uesemo}\\
&&\sup\limits_{t\in[0,T]}|u^{0,\xi}_{u_0}(t)|^2_\mH+\int_0^T\|u^{0,\xi}_{u_0}(t)\|^2_\mV\dif t\leq C(1+|u_0|^2_\mH), \mP-a.s..\label{u0semo}
\ee
\el
\begin{proof}
Since the proofs of (\ref{uesemo}) and (\ref{u0semo}) are similar, we only prove (\ref{uesemo}).

By the It\^o formula, we have that
\ce
|u^{\e,\xi}_{u_0}(t)|^2_\mH&=&|u_0|^2_\mH+2\int_0^t\<u^{\e,\xi}_{u_0}(s),\frac{\partial^2 u^{\e,\xi}_{u_0}(s)}{\partial x^2}\>\dif s+2\int_0^t\<u^{\e,\xi}_{u_0}(s),\frac{\partial g(s,u^{\e,\xi}_{u_0}(s))}{\partial x}\>\dif s\\
&&+2\int_0^t\<u^{\e,\xi}_{u_0}(s),f(s,\cdot,u^{\e,\xi}_{u_0}(s))\>\dif s+2\sqrt\e\sum\limits_{j=1}^d\int_0^t\<u^{\e,\xi}_{u_0}(s),\sigma_j(s,\cdot,u^{\e,\xi}_{u_0}(s))\>\dif W_j(s)\\
&&+2\sum\limits_{j=1}^d\int_0^t\<u^{\e,\xi}_{u_0}(s),\sigma_j(s,\cdot,u^{\e,\xi}_{u_0}(s))\>\xi(s)\dif s+2\int_0^t\int_0^1u^{\e,\xi}_{u_0}(s,x)K(\dif s,\dif x)\\
&&+\e\sum\limits_{j=1}^d\int_0^t|\sigma_j(s,\cdot,u^{\e,\xi}_{u_0}(s))|^2_\mH\dif s.
\de

Next, integration by parts implies that
\ce
2\int_0^t\<u^{\e,\xi}_{u_0}(s),\frac{\partial^2 u^{\e,\xi}_{u_0}(s)}{\partial x^2}\>\dif s=-2\int_0^t\|u^{\e,\xi}_{u_0}(s)\|^2_\mV\dif s.
\de
Set $G(s,z):=\int_0^zg(s,y)\dif y$, and it holds that
\ce
2\int_0^t\<u^{\e,\xi}_{u_0}(s),\frac{\partial g(s,u^{\e,\xi}_{u_0}(s))}{\partial x}\>\dif s&=&-2\int_0^t\<\frac{\partial u^{\e,\xi}_{u_0}(s)}{\partial x}, g(s,u^{\e,\xi}_{u_0}(s))\>\dif s\\
&=&-2\int_0^t\int_0^1\dif G(s,u^{\e,\xi}_{u_0}(s,x))\dif s\\
&=&-2\int_0^t\[G(s,u^{\e,\xi}_{u_0}(s,1))-G(s,u^{\e,\xi}_{u_0}(s,0))\]\dif s\\
&=&0.
\de
$(\mathbf{H}_{f})$ and the H\"older inequality yield that
\ce
2\int_0^t\<u^{\e,\xi}_{u_0}(s),f(s,\cdot,u^{\e,\xi}_{u_0}(s))\>\dif s&\leq& 2\int_0^t\int_0^1|u^{\e,\xi}_{u_0}(s,x)|L^{1/2}_f(1+|u^{\e,\xi}_{u_0}(s,x)|)\dif x\dif s\\
&\leq&4L^{1/2}_f\int_0^t\int_0^1(1+|u^{\e,\xi}_{u_0}(s,x)|^2)\dif x\dif s\\
&\leq&4L^{1/2}_f\int_0^t(1+|u^{\e,\xi}_{u_0}(s)|_\mH^2)\dif s,
\de
where we use the fact that $|z|\leq 1+|z|^2$ for any $z\in\mR$. The Young inequality implies that 
\ce
2\sum\limits_{j=1}^d\int_0^t\<u^{\e,\xi}_{u_0}(s),\sigma_j(s,\cdot,u^{\e,\xi}_{u_0}(s))\>\xi(s)\dif s&\leq& 2\sum\limits_{j=1}^d\int_0^t|u^{\e,\xi}_{u_0}(s)|_\mH|\sigma_j(s,\cdot,u^{\e,\xi}_{u_0}(s))|_\mH |\xi(s)|\dif s\\
&\leq& 2\sup\limits_{s\in[0,t]}|u^{\e,\xi}_{u_0}(s)|_\mH\(\sum\limits_{j=1}^d\int_0^t|\sigma_j(s,\cdot,u^{\e,\xi}_{u_0}(s))|_\mH |\xi(s)|\dif s\)\\
&\leq&\frac{1}{4}\sup\limits_{s\in[0,t]}|u^{\e,\xi}_{u_0}(s)|^2_\mH+4dN\int_0^t|\sigma_j(s,\cdot,u^{\e,\xi}_{u_0}(s))|^2_\mH\dif s,
\de
where in the last inequality we use the fact that $\int_0^T |\xi(s)|^2\dif s\leq N$. Besides, by $(iv)$ in Definition \ref{soludefi}, we know that
$$
2\int_0^t\int_0^1u^{\e,\xi}_{u_0}(s,x)K(\dif s,\dif x)=0.
$$ 

Collecting the above deduction, we have that
\be
&&\frac{3}{4}\sup\limits_{t\in[0,T]}|u^{\e,\xi}_{u_0}(t)|^2_\mH+2\int_0^T\|u^{\e,\xi}_{u_0}(t)\|^2_\mV\dif t\no\\
&\leq&|u_0|^2_\mH+\int_0^T4L^{1/2}_f(1+|u^{\e,\xi}_{u_0}(s)|_\mH^2)\dif s+2\sqrt\e\sum\limits_{j=1}^d\sup\limits_{t\in[0,T]}\left|\int_0^t\<u^{\e,\xi}_{u_0}(s),\sigma_j(s,\cdot,u^{\e,\xi}_{u_0}(s))\>\dif W_j(s)\right|\no\\
&&+(\e+4dN)\sum\limits_{j=1}^d\int_0^T|\sigma_j(s,\cdot,u^{\e,\xi}_{u_0}(s))|^2_\mH\dif s.
\label{sesi}
\ee
The BDG inequality and the Young inequality admit us to obtain that
\ce
&&\frac{3}{4}\mE\sup\limits_{t\in[0,T]}|u^{\e,\xi}_{u_0}(t)|^2_\mH+2\mE\int_0^T\|u^{\e,\xi}_{u_0}(t)\|^2_\mV\dif t\\
&\leq&|u_0|^2_\mH+\int_0^T4L^{1/2}_f(1+\mE|u^{\e,\xi}_{u_0}(s)|_\mH^2)\dif s+2\sqrt\e C\sum\limits_{j=1}^d\mE\(\int_0^T|\<u^{\e,\xi}_{u_0}(s),\sigma_j(s,\cdot,u^{\e,\xi}_{u_0}(s))\>|^2\dif s\)^{1/2}\\
&&+(\e+4dN)\sum\limits_{j=1}^d\mE\int_0^T|\sigma_j(s,\cdot,u^{\e,\xi}_{u_0}(s))|^2_\mH\dif s\\
&\leq&|u_0|^2_\mH+\int_0^T4L^{1/2}_f(1+\mE|u^{\e,\xi}_{u_0}(s)|_\mH^2)\dif s+\frac{1}{4}\mE\left(\sup\limits_{t\in[0,T]}|u^{\e,\xi}_{u_0}(t)|^2_\mH\right)\\
&&+C\sum\limits_{j=1}^d\mE\int_0^T|\sigma_j(s,\cdot,u^{\e,\xi}_{u_0}(s))|^2_\mH\dif s\\
&\leq&|u_0|^2_\mH+\frac{1}{4}\mE\left(\sup\limits_{t\in[0,T]}|u^{\e,\xi}_{u_0}(t)|^2_\mH\right)+C\int_0^T(1+\mE|u^{\e,\xi}_{u_0}(s)|_\mH^2)\dif s.
\de
By the Gronwall inequality, we have (\ref{uesemo}). The proof is complete.
\end{proof}

\bl\label{unummoes}
Under $(\mathbf{H}_{g})$, $(\mathbf{H}_{f})$ and $(\mathbf{H}_{\s})$, for any $R>0, N>0$, there exists a constant $\a_1>0$ such that
\ce
&&\lim\limits_{\e\rightarrow0}\sup\limits_{\|u_0\|_{\mH}\leq R}\sup\limits_{\xi\in\cA^N}\Bigg\{\mE\left[\sup\limits_{t\in[0,T]}\exp\left\{-\a_1\int_0^t(1+\|u^{\e,\xi}_{u_0}(s)\|^2_\mV+\|u^{0,\xi}_{u_0}(s)\|^2_\mV)\dif s\right\}|u^{\e,\xi}_{u_0}(t)-u^{0,\xi}_{u_0}(t)|^2_\mH\right]\\
&&\qquad\quad+\mE\left[\int_0^T\exp\left\{-\a_1\int_0^t(1+\|u^{\e,\xi}_{u_0}(s)\|^2_\mV+\|u^{0,\xi}_{u_0}(s)\|^2_\mV)\dif s\right\}\|u^{\e,\xi}_{u_0}(t)-u^{0,\xi}_{u_0}(t)\|^2_\mV\dif t\right]\Bigg\}=0.
\de
\el
\begin{proof}
First of all, we notice that $u^{\e,\xi}_{u_0}(t)-u^{0,\xi}_{u_0}(t)$ satisfies the following equation in $\mV^*$, $\mP$-a.s.
\ce
u^{\e,\xi}_{u_0}(t)-u^{0,\xi}_{u_0}(t)&=&\int_0^t\frac{\partial^2 \(u^{\e,\xi}_{u_0}(s)-u^{0,\xi}_{u_0}(s)\)}{\partial x^2} \dif s+\int_0^t\frac{\partial \(g(s,u^{\e,\xi}_{u_0}(s))-g(s,u^{0,\xi}_{u_0}(s))\)}{\partial x} \dif s\\
&&+\int_0^t\(f(s,\cdot,u^{\e,\xi}_{u_0}(s))-f(s,\cdot,u^{0,\xi}_{u_0}(s))\)\dif s\\
&&+\sqrt \e\sum\limits_{j=1}^d\int_0^t\sigma_j(s,\cdot,u^{\e,\xi}_{u_0}(s))\dif W_j(s)\\
&&+\sum\limits_{j=1}^d\int_0^t\(\sigma_j(s,\cdot,u^{\e,\xi}_{u_0}(s))-\sigma_j(s,\cdot,u^{0,\xi}_{u_0}(s))\)\xi(s)\dif s\\
&&+\int_0^tK^{\e,\xi}_{u_0}(\dif s,\dif x)-\int_0^tK^{0,\xi}_{u_0}(\dif s,\dif x).
\de
Thus, for any $\a_1>0$, set 
\ce
\Sigma_{\e,0}(t):=\exp\left\{-\a_1\int_0^t(1+\|u^{\e,\xi}_{u_0}(s)\|^2_\mV+\|u^{0,\xi}_{u_0}(s)\|^2_\mV)\dif s\right\},
\de
and by the It\^o formula and integration by parts, it holds that
\be
&&\Sigma_{\e,0}(t)|u^{\e,\xi}_{u_0}(t)-u^{0,\xi}_{u_0}(t)|^2_\mH+2\int_0^t\Sigma_{\e,0}(s)\|u^{\e,\xi}_{u_0}(s)-u^{0,\xi}_{u_0}(s)\|^2_\mV\dif s\no\\
&=&-\a_1\int_0^t\Sigma_{\e,0}(s)(1+\|u^{\e,\xi}_{u_0}(s)\|^2_\mV+\|u^{0,\xi}_{u_0}(s)\|^2_\mV)|u^{\e,\xi}_{u_0}(s)-u^{0,\xi}_{u_0}(s)|^2_\mH\dif s\no\\
&&+2\int_0^t\Sigma_{\e,0}(s)\<u^{\e,\xi}_{u_0}(s)-u^{0,\xi}_{u_0}(s),\frac{\partial \(g(s,u^{\e,\xi}_{u_0}(s))-g(s,u^{0,\xi}_{u_0}(s))\)}{\partial x}\>\dif s\no\\
&&+2\int_0^t\Sigma_{\e,0}(s)\<u^{\e,\xi}_{u_0}(s)-u^{0,\xi}_{u_0}(s),f(s,\cdot,u^{\e,\xi}_{u_0}(s))-f(s,\cdot,u^{0,\xi}_{u_0}(s))\>\dif s\no\\
&&+2\sqrt\e\sum\limits_{j=1}^d\int_0^t\Sigma_{\e,0}(s)\<u^{\e,\xi}_{u_0}(s)-u^{0,\xi}_{u_0}(s),\sigma_j(s,\cdot,u^{\e,\xi}_{u_0}(s))\>\dif W_j(s)\no\\
&&+2\sum\limits_{j=1}^d\int_0^t\Sigma_{\e,0}(s)\<u^{\e,\xi}_{u_0}(s)-u^{0,\xi}_{u_0}(s),\sigma_j(s,\cdot,u^{\e,\xi}_{u_0}(s))-\sigma_j(s,\cdot,u^{0,\xi}_{u_0}(s))\>\xi(s)\dif s\no\\
&&+2\int_0^t\int_0^1\Sigma_{\e,0}(s)\(u^{\e,\xi}_{u_0}(s,x)-u^{0,\xi}_{u_0}(s,x)\)K^{\e,\xi}_{u_0}(\dif s,\dif x)\no\\
&&-2\int_0^t\int_0^1\Sigma_{\e,0}(s)\(u^{\e,\xi}_{u_0}(s,x)-u^{0,\xi}_{u_0}(s,x)\)K^{0,\xi}_{u_0}(\dif s,\dif x)\no\\
&&+\e\sum\limits_{j=1}^d\int_0^t\Sigma_{\e,0}(s)|\sigma_j(s,\cdot,u^{\e,\xi}_{u_0}(s))|^2_\mH\dif s\no\\
&=:&I_1+I_2+I_3+I_4+I_5+I_6+I_7+I_8.
\label{i12345678}
\ee

For $I_2$, by integration by parts, the H\"older inequality and $(\mathbf{H}_{g})$, we have that
\be
I_2&=&-2\int_0^t\Sigma_{\e,0}(s)\<\frac{\partial \(u^{\e,\xi}_{u_0}(s)-u^{0,\xi}_{u_0}(s)\)}{\partial x},g(s,u^{\e,\xi}_{u_0}(s))-g(s,u^{0,\xi}_{u_0}(s))\>\dif s\no\\
&\leq&2\int_0^t\Sigma_{\e,0}(s)\|u^{\e,\xi}_{u_0}(s)-u^{0,\xi}_{u_0}(s)\|_\mV|g(s,u^{\e,\xi}_{u_0}(s))-g(s,u^{0,\xi}_{u_0}(s))|_\mH\dif s\no\\
&\leq&\int_0^t\Sigma_{\e,0}(s)\|u^{\e,\xi}_{u_0}(s)-u^{0,\xi}_{u_0}(s)\|^2_\mV\dif s+\int_0^t\Sigma_{\e,0}(s)|g(s,u^{\e,\xi}_{u_0}(s))-g(s,u^{0,\xi}_{u_0}(s))|^2_\mH\dif s\no\\
&\leq&\int_0^t\Sigma_{\e,0}(s)\|u^{\e,\xi}_{u_0}(s)-u^{0,\xi}_{u_0}(s)\|^2_\mV\dif s\no\\
&&+3L_g^2\int_0^t\Sigma_{\e,0}(s)(1+|u^{\e,\xi}_{u_0}(s)|_{L^\infty([0,1])}^2+|u^{0,\xi}(s)|_{L^\infty([0,1])}^2)|u^{\e,\xi}_{u_0}(s)-u^{0,\xi}_{u_0}(s)|^2_\mH\dif s\no\\
&\leq&\int_0^t\Sigma_{\e,0}(s)\|u^{\e,\xi}_{u_0}(s)-u^{0,\xi}_{u_0}(s)\|^2_\mV\dif s\no\\
&&+3L_g^2\int_0^t\Sigma_{\e,0}(s)(1+\|u^{\e,\xi}_{u_0}(s)\|^2_\mV+\|u^{0,\xi}_{u_0}(s)\|^2_\mV)|u^{\e,\xi}_{u_0}(s)-u^{0,\xi}_{u_0}(s)|^2_\mH\dif s,
\label{i2}
\ee
where the last inequality is based on $|u_{u_0}^{\e,\xi}(s)|_{L^\infty([0,1])}\leq \|u_{u_0}^{\e,\xi}(s)\|_\mV$. For $I_3$, $(\mathbf{H}_{f})$ implies that
\be
I_3\leq 2L_f\int_0^t\Sigma_{\e,0}(s)|u^{\e,\xi}_{u_0}(s)-u^{0,\xi}_{u_0}(s)|^2_\mH\dif s.
\label{i3}
\ee

For $I_5$, by $(\mathbf{H}_{\s})$ and the Young inequality, it holds that
\be
I_5&\leq& 2\sup\limits_{s\in[0,t]}\Sigma^{1/2}_{\e,0}(s)|u^{\e,\xi}_{u_0}(s)-u^{0,\xi}_{u_0}(s)|_\mH\no\\
&&\times\(\sum\limits_{j=1}^d\int_0^t\Sigma^{1/2}_{\e,0}(s)|\sigma_j(s,\cdot,u^{\e,\xi}_{u_0}(s))-\sigma_j(s,\cdot,u^{0,\xi}_{u_0}(s))|_\mH|\xi(s)|\dif s\)\no\\
&\leq&\frac{1}{4}\sup\limits_{s\in[0,t]}\Sigma_{\e,0}(s)|u^{\e,\xi}_{u_0}(s)-u^{0,\xi}_{u_0}(s)|^2_\mH\no\\
&&+4dN\sum\limits_{j=1}^d\int_0^t\Sigma_{\e,0}(s)|\sigma_j(s,\cdot,u^{\e,\xi}_{u_0}(s))-\sigma_j(s,\cdot,u^{0,\xi}_{u_0}(s))|^2_\mH\dif s\no\\
&\leq&\frac{1}{4}\sup\limits_{s\in[0,t]}\Sigma_{\e,0}(s)|u^{\e,\xi}_{u_0}(s)-u^{0,\xi}_{u_0}(s)|^2_\mH\no\\
&&+4dNL_\s\int_0^t\Sigma_{\e,0}(s)|u^{\e,\xi}_{u_0}(s)-u^{0,\xi}_{u_0}(s)|^2_\mH\dif s,
\label{i5}
\ee
where we use the fact that $\int_0^T |\xi(s)|^2\dif s\leq N$.

To deal with $I_6$, we observe that 
\be
I_6&=&2\int_0^t\int_0^1\Sigma_{\e,0}(s)u^{\e,\xi}_{u_0}(s,x)K^{\e,\xi}_{u_0}(\dif s,\dif x)-2\int_0^t\int_0^1\Sigma_{\e,0}(s)u^{0,\xi}_{u_0}(s,x)K^{\e,\xi}_{u_0}(\dif s,\dif x)\no\\
&=&0-2\int_0^t\int_0^1\Sigma_{\e,0}(s)u^{0,\xi}_{u_0}(s,x)K^{\e,\xi}_{u_0}(\dif s,\dif x)\leq 0.
\label{i6}
\ee
By the same deduction to the above inequality, it holds that
\be
I_7\leq 0.
\label{i7}
\ee

For $I_8$, $(\mathbf{H}_{\s})$ yields that
\be
I_8\leq \e L_\s\int_0^t\Sigma_{\e,0}(s)(1+|u^{\e,\xi}_{u_0}(s)|^2_\mH)\dif s.
\label{i8}
\ee

Collecting (\ref{i12345678})-(\ref{i8}) and taking $\a_1\geq 3L_g^2$, we conclude that
\ce
&&\frac{3}{4}\sup\limits_{t\in[0,T]}\Sigma_{\e,0}(t)|u^{\e,\xi}_{u_0}(t)-u^{0,\xi}_{u_0}(t)|^2_\mH+\int_0^T\Sigma_{\e,0}(s)\|u^{\e,\xi}_{u_0}(s)-u^{0,\xi}_{u_0}(s)\|^2_\mV\dif s\\
&\leq&(2L_f+4dNL_\s)\int_0^T\Sigma_{\e,0}(s)|u^{\e,\xi}_{u_0}(s)-u^{0,\xi}_{u_0}(s)|^2_\mH\dif s+\e L_\s\int_0^T\Sigma_{\e,0}(s)(1+|u^{\e,\xi}_{u_0}(s)|^2_\mH)\dif s\\
&&+2\sqrt\e\sum\limits_{j=1}^d\sup\limits_{t\in[0,T]}\left|\int_0^t\Sigma_{\e,0}(s)\<u^{\e,\xi}_{u_0}(s)-u^{0,\xi}_{u_0}(s),\sigma_j(s,\cdot,u^{\e,\xi}_{u_0}(s))\>\dif W_j(s)\right|.
\de
The BDG inequality and the Young inequality imply that
\ce
&&\frac{3}{4}\mE\sup\limits_{t\in[0,T]}\Sigma_{\e,0}(t)|u^{\e,\xi}_{u_0}(t)-u^{0,\xi}_{u_0}(t)|^2_\mH+\mE\int_0^T\Sigma_{\e,0}(s)\|u^{\e,\xi}_{u_0}(s)-u^{0,\xi}_{u_0}(s)\|^2_\mV\dif s\\
&\leq&(2L_f+4dNL_\s)\mE\int_0^T\Sigma_{\e,0}(s)|u^{\e,\xi}_{u_0}(s)-u^{0,\xi}_{u_0}(s)|^2_\mH\dif s+\e L_\s\mE\int_0^T\Sigma_{\e,0}(s)(1+|u^{\e,\xi}_{u_0}(s)|^2_\mH)\dif s\\
&&+2\sqrt\e C\sum\limits_{j=1}^d\mE\left(\int_0^T\Sigma^2_{\e,0}(s)|\<u^{\e,\xi}(s)-u^{0,\xi}(s),\sigma_j(s,\cdot,u^{\e,\xi}_{u_0}(s))\>|^2\dif s\right)^{1/2}\\
&\leq&(2L_f+4dNL_\s)\mE\int_0^T\Sigma_{\e,0}(s)|u^{\e,\xi}_{u_0}(s)-u^{0,\xi}_{u_0}(s)|^2_\mH\dif s+\e L_\s\mE\int_0^T\Sigma_{\e,0}(s)(1+|u^{\e,\xi}_{u_0}(s)|^2_\mH)\dif s\\
&&+\frac{1}{4}\mE\sup\limits_{t\in[0,T]}\Sigma_{\e,0}(t)|u^{\e,\xi}_{u_0}(t)-u^{0,\xi}_{u_0}(t)|^2_\mH+C\e L_\s\mE\int_0^T\Sigma_{\e,0}(s)(1+|u^{\e,\xi}_{u_0}(s)|^2_\mH)\dif s.
\de
By the Gronwall inequality and Lemma \ref{uesefomoes}, it holds that
\ce
&&\mE\sup\limits_{t\in[0,T]}\Sigma_{\e,0}(t)|u^{\e,\xi}_{u_0}(t)-u^{0,\xi}_{u_0}(t)|^2_\mH+\mE\int_0^T\Sigma_{\e,0}(s)\|u^{\e,\xi}_{u_0}(s)-u^{0,\xi}_{u_0}(s)\|^2_\mV\dif s\\
&\leq&C\e L_\s\mE\int_0^T\Sigma_{\e,0}(s)(1+|u^{\e,\xi}_{u_0}(s)|^2_\mH)\dif s\\
&\leq&C\e(1+|u_0|^2_\mH),
\de
where we use the fact that $\Sigma_{\e,0}(s)\leq 1$. So, as $\e\rightarrow 0$, the above inequality implies the required limit. The proof is complete.
\end{proof}

\bl\label{unmcauseq} 
Under $(\mathbf{H}_{g})$, $(\mathbf{H}_{f})$ and $(\mathbf{H}_{\s})$, for any $R>0, N>0$, it holds that for any $\d>0$
\ce
\lim\limits_{\e\rightarrow0}\sup\limits_{\|u_0\|_{\mH}\leq R}\sup\limits_{\xi\in\cA^N}\mP\left\{\sup\limits_{t\in[0,T]}|u^{\e,\xi}_{u_0}(t)-u^{0,\xi}_{u_0}(t)|^2_\mH+\int_0^T\|u^{\e,\xi}_{u_0}(t)-u^{0,\xi}_{u_0}(t)\|^2_\mV\dif t\geq\d\right\}=0.
\de
\el
\begin{proof}
Given $\d>0$, for any $M>0$, it holds that
\ce
&&\mP\left\{\sup\limits_{t\in[0,T]}|u^{\e,\xi}_{u_0}(t)-u^{0,\xi}_{u_0}(t)|^2_\mH+\int_0^T\|u^{\e,\xi}_{u_0}(t)-u^{0,\xi}_{u_0}(t)\|^2_\mV\dif t\geq\d\right\}\\
&=&\mP\bigg\{\sup\limits_{t\in[0,T]}|u^{\e,\xi}_{u_0}(t)-u^{0,\xi}_{u_0}(t)|^2_\mH+\int_0^T\|u^{\e,\xi}_{u_0}(t)-u^{0,\xi}_{u_0}(t)\|^2_\mV\dif t\geq\d,\\
&& \qquad \qquad \a_1\int_0^T(1+\|u^{\e,\xi}_{u_0}(s)\|^2_\mV+\|u^{0,\xi}_{u_0}(s)\|^2_\mV)\dif s<M\bigg\}\\
&&+\mP\bigg\{\sup\limits_{t\in[0,T]}|u^{\e,\xi}_{u_0}(t)-u^{0,\xi}_{u_0}(t)|^2_\mH+\int_0^T\|u^{\e,\xi}_{u_0}(t)-u^{0,\xi}_{u_0}(t)\|^2_\mV\dif t\geq\d,\\
&& \qquad \qquad \a_1\int_0^T(1+\|u^{\e,\xi}_{u_0}(s)\|^2_\mV+\|u^{0,\xi}_{u_0}(s)\|^2_\mV)\dif s\geq M\bigg\}\\
&\leq&\mP\bigg\{\sup\limits_{t\in[0,T]}\Sigma_{\e,0}(t)|u^{\e,\xi}_{u_0}(t)-u^{0,\xi}_{u_0}(t)|^2_\mH+\int_0^T\Sigma_{\e,0}(t)\|u^{\e,\xi}_{u_0}(t)-u^{0,\xi}_{u_0}(t)\|^2_\mV\dif t\geq e^{-M}\d\bigg\}\\
&&+\mP\bigg\{\a_1\int_0^T(1+\|u^{\e,\xi}_{u_0}(s)\|^2_\mV+\|u^{0,\xi}_{u_0}(s)\|^2_\mV)\dif s\geq M\bigg\}\\
&\leq&\frac{e^M}{\d}\left[\mE\sup\limits_{t\in[0,T]}\Sigma_{\e,0}(t)|u^{\e,\xi}_{u_0}(t)-u^{0,\xi}_{u_0}(t)|^2_\mH+\mE\int_0^T\Sigma_{\e,0}(t)\|u^{\e,\xi}_{u_0}(t)-u^{0,\xi}_{u_0}(t)\|^2_\mV\dif t\right]\\
&&+\frac{\a_1}{M}\mE\int_0^T(1+\|u^{\e,\xi}_{u_0}(s)\|^2_\mV+\|u^{0,\xi}_{u_0}(s)\|^2_\mV)\dif s\\
&\leq&\frac{e^M}{\d}\left[\mE\sup\limits_{t\in[0,T]}\Sigma_{\e,0}(t)|u^{\e,\xi}_{u_0}(t)-u^{0,\xi}_{u_0}(t)|^2_\mH+\mE\int_0^T\Sigma_{\e,0}(t)\|u^{\e,\xi}_{u_0}(t)-u^{0,\xi}_{u_0}(t)\|^2_\mV\dif t\right]\\
&&+\frac{\a_1 \(T+2C(1+|u_0|^2_\mH)\)}{M},
\de
where we use (\ref{uesemo}) and (\ref{u0semo}) in the last inequality. As $\e\rightarrow0$ first and then $M\rightarrow\infty$, by Lemma \ref{unummoes} we conclude that
\ce
\lim\limits_{\e\rightarrow0}\sup\limits_{\|u_0\|_{\mH}\leq R}\sup\limits_{\xi\in\cA^N}\mP\left\{\sup\limits_{t\in[0,T]}|u^{\e,\xi}_{u_0}(t)-u^{0,\xi}_{u_0}(t)|^2_\mH+\int_0^T\|u^{\e,\xi}_{u_0}(t)-u^{0,\xi}_{u_0}(t)\|^2_\mV\dif t\geq\d\right\}=0.
\de
This completes the proof of Lemma \ref{unmcauseq}.
\end{proof}

It is the position to prove Theorem \ref{fwuldpsbe}.

{\bf Proof of Theorem \ref{fwuldpsbe}.} By Lemma \ref{unmcauseq}, we know that Condition \ref{cond} holds. Then Theorem \ref{uldpfwjuth} implies Theorem \ref{fwuldpsbe}.

\section{Proof of Theorem \ref{dzuldpsbe}}\label{dzuldpproo}

In this section, we prove Theorem \ref{dzuldpsbe}. In Subsection \ref{compleve} and \ref{contleve}, we justify the compactness of level sets and continuity of level sets in the Hausdorff metric, respectively, which implies that Condition \ref{abcd} holds, and then show Theorem \ref{dzuldpsbe} in terms of Theorem \ref{uldpdzjuth}.

\subsection{Compactness of level sets}\label{compleve}

In this subsection, we prove that for any $u_0\in\mH$ and $M\geq 0$, the level sets $\Phi_{u_0}(M):=\{\phi\in C([0,T],\mH)\cap L^2([0,T],\mV): \Lambda_{u_0}(\phi)\leq M \}$ is compact in $C([0,T],\mH)\cap L^2([0,T],\mV)$.

First of all, for any $n\in\mN$, consider the following penalized Burgers type equation associated with Eq.(\ref{conequ0h}):
\be\left\{\begin{array}{ll}
\dif u^{0,h,n}_{u_0}(t,x)=\frac{\partial^2 u^{0,h,n}_{u_0}(t,x)}{\partial x^2} \dif t+\frac{\partial g(t,u^{0,h,n}_{u_0}(t, x))}{\partial x} \dif t+f(t,x,u^{0,h,n}_{u_0}(t,x))\dif t\\
\qquad\qquad\quad +\sum\limits_{j=1}^d\sigma_j(t,x,u^{0,h,n}_{u_0}(t,x))h_j(t) \dif t+nu^{0,h,n}_{u_0}(t,x)^-\dif t, \quad t\in[0,T], x\in[0,1],\\
u^{0,h,n}_{u_0}(0, \cdot)=u_0(\cdot) \geq 0, \\
u^{0,h,n}_{u_0}(t, 0)=u^{0,h,n}_{u_0}(t, 1)=0,
\end{array}
\right.
\label{psbe}
\ee
where $u^{0,h,n}_{u_0}(s,x)^-:=-\min\{u^{0,h,n}_{u_0}(s,x),0\}$. Under $(\mathbf{H}_{g})$, $(\mathbf{H}_{f})$ and $(\mathbf{H}_{\s})$, by the similar deduction to that for Theorem 2.3 in \cite{gn}, we know that the above equation has a unique solution $u^{0,h,n}_{u_0}$ with
\ce
\sup\limits_{t\in[0,T]}|u^{0,h,n}_{u_0}(t)|^2_\mH+\int_0^T\|u^{0,h,n}_{u_0}(t)\|^2_\mV\dif t<\infty.
\de
Moreover, the following lemma characterizes the relationship between $u^{0,h}_{u_0}$ and $u^{0,h,n}_{u_0}$.

\bl
Under $(\mathbf{H}_{g})$, $(\mathbf{H}_{f})$ and $(\mathbf{H}_{\s})$, for any $N>0$ it holds that
\be
\lim\limits_{n\rightarrow\infty}\sup\limits_{h\in{\bf D}^N}\left\{\sup\limits_{t\in[0,T]}|u^{0,h,n}_{u_0}(t)-u^{0,h}_{u_0}(t)|^2_\mH+\int_0^T\|u^{0,h,n}_{u_0}(t)-u^{0,h}_{u_0}(t)\|^2_\mV\dif t\right\}=0.
\label{0hn0hdiff}
\ee
\el
\begin{proof}
If we prove that
\be
\lim\limits_{n,m\rightarrow\infty}\sup\limits_{h\in{\bf D}^N}\left\{\sup\limits_{t\in[0,T]}|u^{0,h,n}_{u_0}(t)-u^{0,h,m}_{u_0}(t)|^2_\mH+\int_0^T\|u^{0,h,n}_{u_0}(t)-u^{0,h,m}_{u_0}(t)\|^2_\mV\dif t\right\}=0,
\label{0hn0hmdiff}
\ee
the completeness of $C([0,T],\mH)\cap L^2([0,T],\mV)$ and the deduction similar to that in \cite[Theorem 3.1]{zhangt} imply (\ref{0hn0hdiff}). And the proof of (\ref{0hn0hmdiff}) is similar to that for Lemma 3.2 in \cite{zhangt}. Therefore, we omit it. So, the proof is complete.
\end{proof}

Next, we study the relationship between $u^{0,h^k,n}_{u_0}$ and $u^{0,h,n}_{u_0}$  as $h^k\rightarrow h$ in ${\bf D}^N$. For this, we prepare two following lemmas.

\bl
Suppose that $(\mathbf{H}_{g})$, $(\mathbf{H}_{f})$ and $(\mathbf{H}_{\s})$ hold. If $\{h^k, k\in\mN\}\subset{\bf D}^N$, it holds that
\be
\sup\limits_{k}\left\{\sup\limits_{t\in[0,T]}|u^{0,h^k,n}_{u_0}(t)|^2_\mH+\int_0^T\|u^{0,h^k,n}_{u_0}(t)\|^2_\mV\dif t\right\}\leq C(1+|u_0|^2_\mH),
\label{uhkes}
\ee
where $u^{0,h^k,n}_{u_0}$ is a unique solution to Eq.(\ref{psbe}) with replacing $h$ by $h^k$.
\el
\begin{proof}
By the chain rule and integration by parts, we have that
\ce
&&|u^{0,h^k,n}_{u_0}(t)|^2_\mH+2\int_0^t\|u^{0,h^k,n}_{u_0}(s)\|^2_\mV\dif s\\
&=&|u_0|^2_\mH+2\int_0^t\<u^{0,h^k,n}_{u_0}(s),\frac{\partial g(s,u^{0,h^k,n}_{u_0}(s))}{\partial x}\>\dif s+2\int_0^t\<u^{0,h^k,n}_{u_0}(s),f(s,\cdot,u^{0,h^k,n}_{u_0}(s))\>\dif s\\
&&+2\sum\limits_{j=1}^d\int_0^t\<u^{0,h^k,n}_{u_0}(s),\sigma_j(s,\cdot,u^{0,h^k,n}_{u_0}(s))\>h^k_j(s)\dif s+2n\int_0^t\<u^{0,h^k,n}_{u_0}(s),u^{0,h^k,n}_{u_0}(s)^-\>\dif s.
\de

Next, by the similar deduction to that in Lemma \ref{uesefomoes}, it holds that
\ce
&&2\int_0^t\<u^{0,h^k,n}_{u_0}(s),\frac{\partial g(s,u^{0,h^k,n}_{u_0}(s))}{\partial x}\>\dif s=0,\\
&&2\int_0^t\<u^{0,h^k,n}_{u_0}(s),f(s,\cdot,u^{0,h^k,n}_{u_0}(s))\>\dif s\leq 4L^{1/2}_f\int_0^t(1+|u^{0,h^k,n}_{u_0}(s)|_\mH^2)\dif s,\\
&&2\sum\limits_{j=1}^d\int_0^t\<u^{0,h^k,n}_{u_0}(s),\sigma_j(s,\cdot,u^{0,h^k,n}_{u_0}(s))\>h^k_j(s)\dif s\\
&\leq& \frac{1}{2}\sup\limits_{s\in[0,t]}|u^{0,h^k,n}_{u_0}(s)|^2_\mH+2dNL_\s\int_0^t(1+|u^{0,h^k,n}_{u_0}(s)|_\mH^2)\dif s.
\de
Besides, it is easy to see that
\ce
2n\int_0^t\<u^{0,h^k,n}_{u_0}(s),u^{0,h^k,n}_{u_0}(s)^-\>\dif s&=&2n\int_0^t\<u^{0,h^k,n}_{u_0}(s)^+-u^{0,h^k,n}_{u_0}(s)^-,u^{0,h^k,n}_{u_0}(s)^-\>\dif s\\
&=&2n\int_0^t\<-u^{0,h^k,n}_{u_0}(s)^-,u^{0,h^k,n}_{u_0}(s)^-\>\dif s\\
&=&-2n\int_0^t|u^{0,h^k,n}_{u_0}(s)^-|_\mH^2\dif s\leq 0.
\de

Collecting the above deduction, we have that
\ce
&&\frac{1}{2}\sup\limits_{t\in[0,T]}|u^{0,h^k,n}_{u_0}(t)|^2_\mH+2\int_0^T\|u^{0,h^k,n}_{u_0}(t)\|^2_\mV\dif t\no\\
&\leq&|u_0|^2_\mH+(4L^{1/2}_f+2dNL_\s)\int_0^T(1+|u^{0,h^k,n}_{u_0}(s)|_\mH^2)\dif s.
\de
By the Gronwall inequality, we have (\ref{uhkes}). The proof is complete.
\end{proof}

Now, for any $\a_2>0$, set
\ce
&&\Sigma_{n}(t):=\exp\left\{-\a_2\int_0^t(1+\|u^{0,h^k,n}_{u_0}(s)\|^2_\mV+\|u^{0,h,n}_{u_0}(s)\|^2_\mV)\dif s\right\},\\
&&H^k_j(t):=\int_0^t\Sigma_{n}(s)\<u^{0,h^k,n}_{u_0}(s)-u^{0,h,n}_{u_0}(s),\sigma_j(s,\cdot,u^{0,h,n}_{u_0}(s))\>\(h^k_j(s)-h_j(s)\)\dif s, 
\de
where $j=1,2,\dots,d$, and $\{h, h^k; k\in\mN\}\subset{\bf D}^N$.

\bl
Suppose that $(\mathbf{H}_{g})$, $(\mathbf{H}_{f})$ and $(\mathbf{H}_{\s})$ hold. If $h^k\rightarrow h$ in ${\bf D}^N$ as $k\rightarrow \infty$, it holds that
\be
\lim\limits_{k\rightarrow\infty}\sup\limits_{t\in[0,T]}|H^k_j(t)|=0.
\label{limhkj}
\ee
\el
\begin{proof}
First of all, for $0\leq r<t\leq T$, the H\"older inequality implies that
\ce
|H^k_j(t)-H^k_j(r)|&=&\left|\int_r^t\Sigma_{n}(s)\<u^{0,h^k,n}_{u_0}(s)-u^{0,h,n}_{u_0}(s),\sigma_j(s,\cdot,u^{0,h,n}_{u_0}(s))\>\(h^k_j(s)-h_j(s)\)\dif s\right|\\
&\leq&2N^{1/2}\left(\int_r^t\Sigma^2_{n}(s)|u^{0,h^k,n}_{u_0}(s)-u^{0,h,n}_{u_0}(s)|^2_\mH|\sigma_j(s,\cdot,u^{0,h,n}_{u_0}(s))|^2_\mH\dif s\right)^{1/2}\\
&\leq&C\(\sup\limits_{s\in[0,T]}|u^{0,h^k,n}_{u_0}(s)|^2_\mH+\sup\limits_{s\in[0,T]}|u^{0,h,n}_{u_0}(s)|^2_\mH\)^{1/2}\(1+\sup\limits_{s\in[0,T]}|u^{0,h,n}_{u_0}(s)|^2_\mH\)^{1/2}\\
&&\times (t-r)^{1/2}\\
&\leq&C(1+|u_0|^2_\mH)(t-r)^{1/2},
\de
where we use (\ref{uhkes}) and $\Sigma^2_{n}(s)\leq 1$.

Next, let $r=0$, and it holds that
\ce
\sup\limits_{k\in\mN}\sup\limits_{t\in[0,T]}|H^k_j(t)|\leq C(1+|u_0|^2_\mH)T^{1/2}.
\de
Besides, for any $r,t\in[0,T]$
\ce
|H^k_j(t)-H^k_j(r)|\leq C(1+|u_0|^2_\mH)|t-r|^{1/2}.
\de
Combining the above deduction, by the Ascoli-Arzel\'a lemma, we conclude that $\{H^k_j, k\in\mN\}$ is relatively compact in $C([0,T],\mR)$.

Finally, note that
\ce
&&\int_0^T|\Sigma_{n}(s)\<u^{0,h^k,n}_{u_0}(s)-u^{0,h,n}_{u_0}(s),\sigma_j(s,\cdot,u^{0,h,n}_{u_0}(s))\>|^2\dif s\\
&\leq&CT\(\sup\limits_{s\in[0,T]}|u^{0,h^k,n}_{u_0}(s)|^2_\mH+\sup\limits_{s\in[0,T]}|u^{0,h,n}_{u_0}(s)|^2_\mH\)\(1+\sup\limits_{s\in[0,T]}|u^{0,h,n}_{u_0}(s)|^2_\mH\)\\
&\leq&CT(1+|u_0|^2_\mH)^2<\infty.
\de
Thus, since $h^k\rightarrow h$ in ${\bf D}^N$ as $k\rightarrow \infty$, 
$$
\lim\limits_{k\rightarrow\infty}H^k_j(t)=0,
$$
which together with the relative compactness of $\{H^k_j, k\in\mN\}$ yields that
\ce
\lim\limits_{k\rightarrow\infty}\sup\limits_{t\in[0,T]}|H^k_j(t)|=0.
\de
The proof is complete.
\end{proof}

Here we prove the convergence of $u^{0,h^k,n}_{u_0}$ to $u^{0,h,n}_{u_0}$ as $h^k\rightarrow h$ in ${\bf D}^N$. 

\bl
Suppose that $(\mathbf{H}_{g})$, $(\mathbf{H}_{f})$ and $(\mathbf{H}_{\s})$ hold. If $h^k\rightarrow h$ in ${\bf D}^N$ as $k\rightarrow \infty$, we have that
\be
\lim\limits_{k\rightarrow\infty}\left\{\sup\limits_{t\in[0,T]}|u^{0,h^k,n}_{u_0}(t)-u^{0,h,n}_{u_0}(t)|^2_\mH+\int_0^T\|u^{0,h^k,n}_{u_0}(t)-u^{0,h,n}_{u_0}(t)\|^2_\mV\dif t\right\}=0.
\label{0hkn0hndiff}
\ee
\el
\begin{proof}
By the chain rule and integration by parts, it holds that
\be
&&\Sigma_{n}(t)|u^{0,h^k,n}_{u_0}(t)-u^{0,h,n}_{u_0}(t)|^2_\mH+2\int_0^t\Sigma_{n}(s)\|u^{0,h^k,n}_{u_0}(s)-u^{0,h,n}_{u_0}(s)\|^2_\mV\dif s\no\\
&=&-\a_2\int_0^t\Sigma_{n}(s)(1+\|u^{0,h^k,n}_{u_0}(s)\|^2_\mV+\|u^{0,h,n}_{u_0}(s)\|^2_\mV)|u^{0,h^k,n}_{u_0}(s)-u^{0,h,n}_{u_0}(s)|^2_\mH\dif s\no\\
&&+2\int_0^t\Sigma_{n}(s)\<u^{0,h^k,n}_{u_0}(s)-u^{0,h,n}_{u_0}(s),\frac{\partial \(g(s,u^{0,h^k,n}_{u_0}(s))-g(s,u^{0,h,n}_{u_0}(s))\)}{\partial x}\>\dif s\no\\
&&+2\int_0^t\Sigma_{n}(s)\<u^{0,h^k,n}_{u_0}(s)-u^{0,h,n}_{u_0}(s),f(s,\cdot,u^{0,h^k,n}_{u_0}(s))-f(s,\cdot,u^{0,h,n}_{u_0}(s))\>\dif s\no\\
&&+2\sum\limits_{j=1}^d\int_0^t\Sigma_{n}(s)\<u^{0,h^k,n}_{u_0}(s)-u^{0,h,n}_{u_0}(s),\sigma_j(s,\cdot,u^{0,h^k,n}_{u_0}(s))-\sigma_j(s,\cdot,u^{0,h,n}_{u_0}(s))\>h^k_j(s)\dif s\no\\
&&+2\sum\limits_{j=1}^d\int_0^t\Sigma_{n}(s)\<u^{0,h^k,n}_{u_0}(s)-u^{0,h,n}_{u_0}(s),\sigma_j(s,\cdot,u^{0,h,n}_{u_0}(s))\>\(h^k_j(s)-h_j(s)\)\dif s\no\\
&&+2\int_0^t\Sigma_{n}(s)\<u^{0,h^k,n}_{u_0}(s)-u^{0,h,n}_{u_0}(s),nu^{0,h^k,n}_{u_0}(s)^--nu^{0,h,n}_{u_0}(s)^-\>\dif s\no\\
&=:&J_1+J_2+J_3+J_4+J_5+J_6.
\label{j123456}
\ee

For $J_2$, $J_3$ and $J_4$, by the similar deduction to that for (\ref{i2}), (\ref{i3}) and (\ref{i5}), we have that
\be
J_2&\leq&\int_0^t\Sigma_{n}(s)\|u^{0,h^k,n}_{u_0}(s)-u^{0,h,n}_{u_0}(s)\|^2_\mV\dif s\no\\
&&+3L_g^2\int_0^t\Sigma_{n}(s)(1+\|u^{0,h^k,n}_{u_0}(s)\|^2_\mV+\|u^{0,h,n}_{u_0}(s)\|^2_\mV)|u^{0,h^k,n}_{u_0}(s)-u^{0,h,n}_{u_0}(s)|^2_\mH\dif s,\label{j2}\\
J_3&\leq& 2L_f\int_0^t\Sigma_{n}(s)|u^{0,h^k,n}_{u_0}(s)-u^{0,h,n}_{u_0}(s)|^2_\mH\dif s,\label{j3}\\
J_4&\leq& \frac{1}{2}\sup\limits_{s\in[0,t]}\Sigma_{n}(s)|u^{0,h^k,n}_{u_0}(s)-u^{0,h,n}_{u_0}(s)|^2_\mH\no\\
&&+2dN L_\s\int_0^t\Sigma_{n}(s)|u^{0,h^k,n}_{u_0}(s)-u^{0,h,n}_{u_0}(s)|^2_\mH\dif s.
\label{j4}
\ee

To deal with $J_6$, we observe that 
\ce
|u^{0,h^k,n}_{u_0}(s,x)^--u^{0,h,n}_{u_0}(s,x)^-|\leq |u^{0,h^k,n}_{u_0}(s,x)-u^{0,h,n}_{u_0}(s,x)|,
\de
which implies that
\be
J_6\leq 2n\int_0^t\Sigma_{n}(s)|u^{0,h^k,n}_{u_0}(s)-u^{0,h,n}_{u_0}(s)|^2_\mH\dif s.
\label{j6}
\ee

Collecting (\ref{j123456})-(\ref{j6}) and letting $\a_2\geq 3L_g^2$, we conclude that
\ce
&&\frac{1}{2}\sup\limits_{t\in[0,T]}\Sigma_{n}(t)|u^{0,h^k,n}_{u_0}(t)-u^{0,h,n}_{u_0}(t)|^2_\mH+\int_0^T\Sigma_{n}(t)\|u^{0,h^k,n}_{u_0}(t)-u^{0,h,n}_{u_0}(t)\|^2_\mV\dif t\\
&\leq&(2L_f+2dN L_\s+2n)\int_0^T\Sigma_{n}(s)|u^{0,h^k,n}_{u_0}(s)-u^{0,h,n}_{u_0}(s)|^2_\mH\dif s\\
&&+2\sum\limits_{j=1}^d\sup\limits_{t\in[0,T]}\left|\int_0^t\Sigma_{n}(s)\<u^{0,h^k,n}_{u_0}(s)-u^{0,h,n}_{u_0}(s),\sigma_j(s,\cdot,u^{0,h,n}_{u_0}(s))\>\(h^k_j(s)-h_j(s)\)\dif s\right|.
\de
By the Gronwall inequality, it holds that
\ce
&&\sup\limits_{t\in[0,T]}\Sigma_{n}(t)|u^{0,h^k,n}_{u_0}(t)-u^{0,h,n}_{u_0}(t)|^2_\mH+\int_0^T\Sigma_{n}(t)\|u^{0,h^k,n}_{u_0}(t)-u^{0,h,n}_{u_0}(t)\|^2_\mV\dif t\no\\
&\leq&C\sum\limits_{j=1}^d\sup\limits_{t\in[0,T]}\left|\int_0^t\Sigma_{n}(s)\<u^{0,h^k,n}_{u_0}(s)-u^{0,h,n}_{u_0}(s),\sigma_j(s,\cdot,u^{0,h,n}_{u_0}(s))\>\(h^k_j(s)-h_j(s)\)\dif s\right|.
\de
Taking the limit on two sides of the above inequality, by (\ref{limhkj}) we obtain that
\be
\lim\limits_{k\rightarrow\infty}\left\{\sup\limits_{t\in[0,T]}\Sigma_{n}(t)|u^{0,h^k,n}_{u_0}(t)-u^{0,h,n}_{u_0}(t)|^2_\mH+\int_0^T\Sigma_{n}(t)\|u^{0,h^k,n}_{u_0}(t)-u^{0,h,n}_{u_0}(t)\|^2_\mV\dif t\right\}=0.
\label{sighhknhndiff}
\ee

Finally, note that 
\ce
\sup\limits_{t\in[0,T]}\Sigma^{-1}_{n}(t)&=&\exp\left\{\a_2\int_0^T(1+\|u^{0,h^k,n}_{u_0}(s)\|^2_\mV+\|u^{0,h,n}_{u_0}(s)\|^2_\mV)\dif s\right\}\\
&\leq&\exp\left\{\a_2(T+C(1+|u_0|^2_\mH))\right\}.
\de
Thus, it holds that
\ce
&&\sup\limits_{t\in[0,T]}|u^{0,h^k,n}_{u_0}(t)-u^{0,h,n}_{u_0}(t)|^2_\mH+\int_0^T\|u^{0,h^k,n}_{u_0}(t)-u^{0,h,n}_{u_0}(t)\|^2_\mV\dif t\\
&\leq&\(\sup\limits_{t\in[0,T]}\Sigma_{n}(t)|u^{0,h^k,n}_{u_0}(t)-u^{0,h,n}_{u_0}(t)|^2_\mH\)\sup\limits_{t\in[0,T]}\Sigma^{-1}_{n}(t)\\
&&+\(\int_0^T\Sigma_{n}(t)\|u^{0,h^k,n}_{u_0}(t)-u^{0,h,n}_{u_0}(t)\|^2_\mV\dif t\)\sup\limits_{t\in[0,T]}\Sigma^{-1}_{n}(t)\\
&\leq&\left\{\sup\limits_{t\in[0,T]}\Sigma_{n}(t)|u^{0,h^k,n}_{u_0}(t)-u^{0,h,n}_{u_0}(t)|^2_\mH+\int_0^T\Sigma_{n}(t)\|u^{0,h^k,n}_{u_0}(t)-u^{0,h,n}_{u_0}(t)\|^2_\mV\dif t\right\}\\
&&\cdot \exp\left\{\a_2(T+C(1+|u_0|^2_\mH))\right\},
\de
which together with (\ref{sighhknhndiff}) implies the required result. 
\end{proof}

Now, keeping (\ref{0hn0hdiff}) and (\ref{0hkn0hndiff}) in mind, we obtain the following result.

\bp\label{rategood}
Suppose that $(\mathbf{H}_{g})$, $(\mathbf{H}_{f})$ and $(\mathbf{H}_{\s})$ hold. Then for any $u_0\in\mH$ and $M\geq 0$, $\Phi_{u_0}(M)$ is compact in $C([0,T],\mH)\cap L^2([0,T],\mV)$.
\ep
\begin{proof}
First of all, we mention that 
\ce
\Phi_{u_0}(M)&=&\{\phi\in C([0,T],\mH)\cap L^2([0,T],\mV): \Lambda_{u_0}(\phi)\leq M\}\\
&=&\{u^{0,h}_{u_0}\in C([0,T],\mH)\cap L^2([0,T],\mV): h\in {\bf D}^{2M}\}.
\de
Moreover, ${\bf D}^{2M}$ is compact under the weak topology in $L^2\left([0, T], \mR^d\right)$. Thus, if we prove that the mapping ${\bf D}^{2M}\ni h\mapsto u^{0,h}_{u_0}\in C([0,T],\mH)\cap L^2([0,T],\mV)$ is continuous, $\Phi_{u_0}(M)$ is compact in $C([0,T],\mH)\cap L^2([0,T],\mV)$. 

In the following, we prove that if $h^k\rightarrow h$ in ${\bf D}^{2M}$ as $k\rightarrow \infty$, $u^{0,h^k}_{u_0}\rightarrow u^{0,h}_{u_0}$ in $C([0,T],\mH)\cap L^2([0,T],\mV)$. In fact, for $u^{0,h^k}_{u_0}, u^{0,h}_{u_0}\in C([0,T],\mH)$,
\ce
&&\sup\limits_{t\in[0,T]}|u^{0,h^k}_{u_0}(t)-u^{0,h}_{u_0}(t)|^2_\mH\\
&\leq& 3\sup\limits_{t\in[0,T]}|u^{0,h^k}_{u_0}(t)-u^{0,h^k,n}_{u_0}(t)|^2_\mH+3\sup\limits_{t\in[0,T]}|u^{0,h^k,n}_{u_0}(t)-u^{0,h,n}_{u_0}(t)|^2_\mH\\
&&+3\sup\limits_{t\in[0,T]}|u^{0,h,n}_{u_0}(t)-u^{0,h}_{u_0}(t)|^2_\mH\\
&\leq& 6\sup\limits_{h\in {\bf D}^{2M}}\sup\limits_{t\in[0,T]}|u^{0,h}_{u_0}(t)-u^{0,h,n}_{u_0}(t)|^2_\mH+3\sup\limits_{t\in[0,T]}|u^{0,h^k,n}_{u_0}(t)-u^{0,h,n}_{u_0}(t)|^2_\mH.
\de   
Letting $k\rightarrow\infty$ first and then $n\rightarrow\infty$, (\ref{0hn0hdiff}) and (\ref{0hkn0hndiff}) imply that $u^{0,h^k}_{u_0}\rightarrow u^{0,h}_{u_0}$ in $C([0,T],\mH)$. By the same deduction to the above, we can obtain that $u^{0,h^k}_{u_0}\rightarrow u^{0,h}_{u_0}$ in $L^2([0,T],\mV)$. The proof is complete.
\end{proof}

\subsection{Continuity of the level sets in the Hausdorff metric}\label{contleve}

In this subsection, we prove that the level sets are continuous in the Hausdorff metric. First we prepare the following result.

\bl\label{u0nu0limi}
Suppose that $(\mathbf{H}_{g})$, $(\mathbf{H}_{f})$ and $(\mathbf{H}_{\s})$ hold. If for any $N\in\mN$ and $h\in {\bf D}^N$, $u_{0,n}\rightarrow u_0$ in $\mH$ as $n\rightarrow\infty$, it holds that
\ce
\lim\limits_{n\rightarrow\infty}\left\{\sup\limits_{t\in[0,T]}|u^{0,h}_{u_{0,n}}(t)-u^{0,h}_{u_0}(t)|^2_\mH+\int_0^T\|u^{0,h}_{u_{0,n}}(t)-u^{0,h}_{u_0}(t)\|^2_\mV\dif t\right\}=0,
\de
where $(u^{0,h}_{u_{0,n}}, K^{0,h}_{u_{0,n}})$ is a unique solution of Eq.(\ref{conequ0h}) with replacing $u_0$ by $u_{0,n}$.
\el
\begin{proof}
First of all, by (\ref{u0semo}), it holds that 
\be
&&\sup\limits_{t\in[0,T]}|u^{0,h}_{u_{0,n}}(t)|^2_\mH+\int_0^T\|u^{0,h}_{u_{0,n}}(t)\|^2_\mV\dif t\leq C(1+|u_{0,n}|^2_\mH),\label{u0nhes}\\
&&\sup\limits_{t\in[0,T]}|u^{0,h}_{u_{0}}(t)|^2_\mH+\int_0^T\|u^{0,h}_{u_{0}}(t)\|^2_\mV\dif t\leq C(1+|u_{0}|^2_\mH).\label{u0hes}
\ee
Then, for any $\a_3>0$, set 
\ce
\Sigma_{0}(t):=\exp\left\{-\a_3\int_0^t(1+\|u^{0,h}_{u_{0,n}}(s)\|^2_\mV+\|u^{0,h}_{u_{0}}(s)\|^2_\mV)\dif s\right\},
\de
and by the chain rule and integration by parts, it holds that
\be
&&\Sigma_{0}(t)|u^{0,h}_{u_{0,n}}(t)-u^{0,h}_{u_0}(t)|^2_\mH+2\int_0^t\Sigma_{0}(s)\|u^{0,h}_{u_{0,n}}(s)-u^{0,h}_{u_0}(s)\|^2_\mV\dif s\no\\
&=&|u_{0,n}-u_0|^2_\mH-\a_3\int_0^t\Sigma_{0}(s)(1+\|u^{0,h}_{u_{0,n}}(s)\|^2_\mV+\|u^{0,h}_{u_{0}}(s)\|^2_\mV)|u^{0,h}_{u_{0,n}}(s)-u^{0,h}_{u_0}(s)|^2_\mH\dif s\no\\
&&+2\int_0^t\Sigma_{0}(s)\<u^{0,h}_{u_{0,n}}(s)-u^{0,h}_{u_0}(s),\frac{\partial \(g(s,u^{0,h}_{u_{0,n}}(s))-g(s,u^{0,h}_{u_{0}}(s))\)}{\partial x}\>\dif s\no\\
&&+2\int_0^t\Sigma_{0}(s)\<u^{0,h}_{u_{0,n}}(s)-u^{0,h}_{u_0}(s),f(s,\cdot, u^{0,h}_{u_{0,n}}(s))-f(s,\cdot, u^{0,h}_{u_{0}}(s))\>\dif s\no\\
&&+2\sum\limits_{j=1}^d\int_0^t\Sigma_{0}(s)\<u^{0,h}_{u_{0,n}}(s)-u^{0,h}_{u_0}(s),\s_j(s,\cdot, u^{0,h}_{u_{0,n}}(s))-\s_j(s,\cdot, u^{0,h}_{u_{0}}(s))\>h_j(s)\dif s\no\\
&&+2\int_0^t\int_0^1\Sigma_{0}(s)\(u^{0,h}_{u_{0,n}}(s)-u^{0,h}_{u_0}(s)\)\dif K^{0,h}_{u_{0,n}}(\dif s, \dif x)\no\\
&&-2\int_0^t\int_0^1\Sigma_{0}(s)\(u^{0,h}_{u_{0,n}}(s)-u^{0,h}_{u_0}(s)\)\dif K^{0,h}_{u_{0}}(\dif s, \dif x)\no\\
&=:&|u_{0,n}-u_0|^2_\mH+\cK_1+\cK_2+\cK_3+\cK_4+\cK_5+\cK_6.
\label{k123456}
\ee
For $\cK_2$, $\cK_3$, $\cK_4$, $\cK_5$ and $J_6$, by the similar deduction to that for (\ref{i2}), (\ref{i3}), (\ref{i5}), (\ref{i6}) and (\ref{i7}), we have that
\be
\cK_2&\leq&\int_0^t\Sigma_{0}(s)\|u^{0,h}_{u_{0,n}}(s)-u^{0,h}_{u_0}(s)\|^2_\mV\dif s\no\\
&&+3L_g^2\int_0^t\Sigma_{0}(s)(1+\|u^{0,h}_{u_{0,n}}(s)\|^2_\mV+\|u^{0,h}_{u_{0}}(s)\|^2_\mV)|u^{0,h}_{u_{0,n}}(s)-u^{0,h}_{u_0}(s)|^2_\mH\dif s,\label{k2}\\
\cK_3&\leq& 2L_f\int_0^t\Sigma_{0}(s)|u^{0,h}_{u_{0,n}}(s)-u^{0,h}_{u_0}(s)|^2_\mH\dif s,\label{k3}\\
\cK_4&\leq& \frac{1}{2}\sup\limits_{s\in[0,t]}\Sigma_{0}(s)|u^{0,h}_{u_{0,n}}(s)-u^{0,h}_{u_0}(s)|^2_\mH+2dN L_\s\int_0^t\Sigma_{0}(s)|u^{0,h}_{u_{0,n}}(s)-u^{0,h}_{u_0}(s)|^2_\mH\dif s,\label{k4}\\
\cK_5&\leq& 0,\label{k5}\\
\cK_6&\leq& 0.\label{k6}
\ee
Combining (\ref{k2})-(\ref{k6}) with (\ref{k123456}), for $\a_3\geq 3L_g^2$, one can infer that
\ce
&&\frac{1}{2}\sup\limits_{t\in[0,T]}\Sigma_{0}(t)|u^{0,h}_{u_{0,n}}(t)-u^{0,h}_{u_0}(t)|^2_\mH+\int_0^T\Sigma_{0}(s)\|u^{0,h}_{u_{0,n}}(s)-u^{0,h}_{u_0}(s)\|^2_\mV\dif s\\
&\leq&|u_{0,n}-u_0|^2_\mH+(2L_f+2dN L_\s)\int_0^T\Sigma_{0}(s)|u^{0,h}_{u_{0,n}}(s)-u^{0,h}_{u_0}(s)|^2_\mH\dif s.
\de
The Gronwall inequality implies that
\ce
\sup\limits_{t\in[0,T]}\Sigma_{0}(t)|u^{0,h}_{u_{0,n}}(t)-u^{0,h}_{u_0}(t)|^2_\mH+\int_0^T\Sigma_{0}(s)\|u^{0,h}_{u_{0,n}}(s)-u^{0,h}_{u_0}(s)\|^2_\mV\dif s\leq C|u_{0,n}-u_0|^2_\mH.
\de

At last, we know that 
\ce
&&\sup\limits_{t\in[0,T]}|u^{0,h}_{u_{0,n}}(t)-u^{0,h}_{u_0}(t)|^2_\mH+\int_0^T\|u^{0,h}_{u_{0,n}}(s)-u^{0,h}_{u_0}(s)\|^2_\mV\dif s\\
&\leq&\(\sup\limits_{t\in[0,T]}\Sigma_{0}(t)|u^{0,h}_{u_{0,n}}(t)-u^{0,h}_{u_0}(t)|^2_\mH\)\sup\limits_{t\in[0,T]}\Sigma^{-1}_{0}(t)\\
&&+\left(\int_0^T\Sigma_{0}(s)\|u^{0,h}_{u_{0,n}}(s)-u^{0,h}_{u_0}(s)\|^2_\mV\dif s\right)\sup\limits_{t\in[0,T]}\Sigma^{-1}_{0}(t)\\
&\leq&C|u_{0,n}-u_0|^2_\mH \sup\limits_{t\in[0,T]}\Sigma^{-1}_{0}(t).
\de
Besides, (\ref{u0nhes}) and (\ref{u0hes}) assure that
\ce
\sup\limits_{t\in[0,T]}\Sigma^{-1}_{0}(t)\leq \exp\left\{\a_3\(T+C(1+|u_{0,n}|^2_\mH)+C(1+|u_{0}|^2_\mH)\)\right\}.
\de
Collecting the above deduction and noticing that $u_{0,n}\rightarrow u_0$ in $\mH$ as $n\rightarrow\infty$, we conclude that
\ce
\lim\limits_{n\rightarrow\infty}\left\{\sup\limits_{t\in[0,T]}|u^{0,h}_{u_{0,n}}(t)-u^{0,h}_{u_0}(t)|^2_\mH+\int_0^T\|u^{0,h}_{u_{0,n}}(t)-u^{0,h}_{u_0}(t)\|^2_\mV\dif t\right\}=0.
\de
The proof is complete.
\end{proof}

\bp\label{contleveha}
Assume that $(\mathbf{H}_{g})$, $(\mathbf{H}_{f})$ and $(\mathbf{H}_{\s})$ hold and $u_{0,n}\rightarrow u_0$ in $\mH$ as $n\rightarrow\infty$. Then for any $M\geq 0$
\ce
&& \lim _{n \rightarrow+\infty} \max\Bigg\{\sup\limits_{\phi\in\Phi_{u_0}(M)}\rho_{C([0,T],\mH)\cap L^2([0,T],\mV)}\(\phi,\Phi_{u_{0,n}}(M)\), \\
 &&\qquad\qquad \sup\limits_{\phi\in\Phi_{u_{0,n}}(M)}\rho_{C([0,T],\mH)\cap L^2([0,T],\mV)}\(\phi,\Phi_{u_{0}}(M)\)\Bigg\}=0.
\de
\ep
\begin{proof}
First of all, if $\phi\in\Phi_{u_0}(M)$, $\Lambda_{u_0}(\phi)\leq M$. By the definition of $\Lambda_{u_0}(\phi)$, there exists a $h\in L^2([0,T],\mR^d)$ such that $\phi=u^{0,h}_{u_0}$ and 
\ce
\frac{1}{2}\int_0^T|h(t)|^2\dif t\leq M,
\de
which also implies that $\Lambda_{u_{0,n}}(u^{0,h}_{u_{0,n}})\leq M$ and $u^{0,h}_{u_{0,n}}\in \Phi_{u_{0,n}}(M)$. So,
\ce
\rho_{C([0,T],\mH)\cap L^2([0,T],\mV)}(\phi, \Phi_{u_{0,n}}(M))&\leq& \rho_{C([0,T],\mH)\cap L^2([0,T],\mV)}(\phi, u^{0,h}_{u_{0,n}})\\
&=&\sup\limits_{t\in[0,T]}|u^{0,h}_{u_{0,n}}(t)-u^{0,h}_{u_0}(t)|^2_\mH+\int_0^T\|u^{0,h}_{u_{0,n}}(t)-u^{0,h}_{u_0}(t)\|^2_\mV\dif t.
\de
Lemma \ref{u0nu0limi} implies that
\ce
\lim\limits_{n\rightarrow\infty}\sup\limits_{\phi\in\Phi_{u_0}(M)}\rho_{C([0,T],\mH)\cap L^2([0,T],\mV)}(\phi, \Phi_{u_{0,n}}(M))=0.
\de

Similarly, it holds that
\ce
\lim\limits_{n\rightarrow\infty}\sup\limits_{\phi\in\Phi_{u_{0,n}}(M)}\rho_{C([0,T],\mH)\cap L^2([0,T],\mV)}(\phi, \Phi_{u_0}(M))=0.
\de
Finally, the above deduction implies the required result.
\end{proof}

It is the position to prove Theorem \ref{dzuldpsbe}.

{\bf Proof of Theorem \ref{dzuldpsbe}.} By Proposition \ref{rategood} and \ref{contleveha}, we conclude that $(c)$ and $(d)$ of Condition \ref{abcd} hold. Then by taking $\mX_0=\mV$ and letting $\mathscr{A}$ be the collection of compact subsets of $\mV$, it is easy to see that $(a)$ and $(b)$ of Condition \ref{abcd} are also right. Thus, by Theorem \ref{uldpdzjuth}, we obtain Theorem \ref{dzuldpsbe}.

\section{Proof of Theorem \ref{averprintheo}}\label{averprintheoproo}

In this section, we prove Theorem \ref{averprintheo}. We first prepare some estimates for penalized stochastic Burgers type equations and show Theorem \ref{averprintheo} through the convergence for penalized stochastic Burgers type equations.

First of all, for any $n\in\mN$, consider the following penalized stochastic Burgers type equation associated with Eq.(\ref{sbe3}):
\be\left\{\begin{array}{ll}
\dif \bar u^{\e,n}(t, x)=\frac{\partial^2 \bar u^{\e,n}(t, x)}{\partial x^2} \dif t+\frac{\partial g(t,\bar u^{\e,n}(t, x))}{\partial x} \dif t+f(\frac{t}{\e},x,\bar u^{\e,n}(t, x))\dif t \\
\qquad\qquad\quad+\sum\limits_{j=1}^d \sigma_j(\frac{t}{\e},x,\bar u^{\e,n}(t, x)) \dif W_j(t)+n\bar u^{\e,n}(t, x)^-\dif t, \quad t\in\mR_+, \quad x \in[0,1],\\
\bar u^{\e,n}(t, x)\geq 0,\\
\bar u^{\e,n}(0, \cdot)=u_0(\cdot) \geq 0, \quad \bar u^{\e,n}(t, 0)=\bar u^{\e,n}(t, 1)=0.
\end{array}
\right.
\label{sbe3pena}
\ee
Then by the similar deduction to that for Theorem 2.3 in \cite{gn}, under $(\mathbf{H}_{g})$, $(\mathbf{H}_{f})$ and $(\mathbf{H}_{\s})$ we know that for $u_0\in\mV$, Eq.(\ref{sbe3pena}) has a unique solution $\bar u^{\e,n}_{u_0}$ in $C([0,T],\mH)\cap L^2([0,T],\mV)$ for any $T>0$. Besides, we present the penalized stochastic Burgers type equation associated with Eq.(\ref{sbe4}):
\be\left\{\begin{array}{ll}
\dif \bar u^{0,n}(t, x)=\frac{\partial^2 \bar u^{0,n}(t, x)}{\partial x^2} \dif t+\frac{\partial g(t,\bar u^{0,n}(t, x))}{\partial x} \dif t+\bar f(x,\bar u^{0,n}(t, x))\dif t \\
\qquad\qquad\quad+\sum\limits_{j=1}^d \bar\sigma_j(x,\bar u^{0,n}(t, x)) \dif W_j(t)+n\bar u^{0,n}(t, x)^-\dif t, \quad t\in\mR_+, \quad x \in[0,1],\\
\bar u^{0,n}(t, x)\geq 0,\\
\bar u^{0,n}(0, \cdot)=u_0(\cdot) \geq 0, \quad \bar u^{0,n}(t, 0)=\bar u^{0,n}(t, 1)=0.
\end{array}
\right.
\label{sbe4pena}
\ee
So, Theorem 2.3 in \cite{gn} admits us to get that for $u_0\in\mV$, Eq.(\ref{sbe4pena}) has a unique solution $\bar u^{0,n}_{u_0}$ under $(\mathbf{H}_{g})$, $(\mathbf{H}_{f})$, $(\mathbf{H}_{\s})$ and $(\mathbf{H}^L_{f,\s})$. We collect some results about $\bar u^{\e,n}_{u_0}$ and $\bar u^{0,n}_{u_0}$ in the following lemma.

\bl\label{baruenbaru0nresu}
Assume that $(\mathbf{H}_{g})$, $(\mathbf{H}_{f})$, $(\mathbf{H}_{\s})$ and $(\mathbf{H}^L_{f,\s})$ hold and $u_0\in\mV$. Then it holds that 

$(i)$ For any $T>0$
\be
&&\sup\limits_{\e\in(0,1),n\in\mN}\mE\left[\sup\limits_{t\in[0,T]}|\bar u^{\e,n}_{u_0}(t)|^2_\mH+\int_0^T\|\bar u^{\e,n}_{u_0}(t)\|^2_\mV\dif t\right]\leq C(1+|u_0|^2_\mH),\label{baruenes2}\\
&&\sup\limits_{\e\in(0,1),n\in\mN}\mE\left[\sup\limits_{t\in[0,T]}|\bar u^{\e,n}_{u_0}(t)|^4_\mH+\left(\int_0^T\|\bar u^{\e,n}_{u_0}(t)\|^2_\mV\dif t\right)^2\right]\leq C(1+|u_0|^4_\mH).\label{baruenes4}
\ee

$(ii)$ For any $T>0$
\be
&&\sup\limits_{n\in\mN}\mE\left[\sup\limits_{t\in[0,T]}|\bar u^{0,n}_{u_0}(t)|^2_\mH+\int_0^T\|\bar u^{0,n}_{u_0}(t)\|^2_\mV\dif t\right]\leq C(1+|u_0|^2_\mH),\label{baru0nes2}\\
&&\sup\limits_{n\in\mN}\mE\left[\sup\limits_{t\in[0,T]}|\bar u^{0,n}_{u_0}(t)|^4_\mH+\left(\int_0^T\|\bar u^{0,n}_{u_0}(t)\|^2_\mV\dif t\right)^2\right]\leq C(1+|u_0|^4_\mH).\label{baru0nes4}
\ee

$(iii)$ For any $T>0$, $\bar u^{\e,n}_{u_0}$ converges to $\bar u^{\e}_{u_0}$ in probability in $C([0,T],\mH)\cap L^2([0,T],\mV)$ as $n\rightarrow \infty$.

$(iv)$ For any $T>0$, $\bar u^{0,n}_{u_0}$ converges to $\bar u^{0}_{u_0}$ in probability in $C([0,T],\mH)\cap L^2([0,T],\mV)$ as $n\rightarrow \infty$.
\el
\begin{proof}
First of all, two estimates (\ref{baruenes2}) and (\ref{baru0nes2}) can be proved by the same deduction to that for (\ref{uesemo}). And by the similar computation to that for \cite[(2.4)]{zhangt}, one can show (\ref{baruenes4}) and (\ref{baru0nes4}). Thus, we have $(i)$ and $(ii)$.

Next, along the same line to one for Lemma 3.2 in \cite{zhangt}, we conclude $(iii)$ and $(iv)$. The proof is complete.
\end{proof}

\bl
Assume that $(\mathbf{H}_{g})$, $(\mathbf{H}_{f})$, $(\mathbf{H}_{\s})$, $(\mathbf{H}^C_{g})$ and $(\mathbf{H}^L_{f,\s})$ hold and $u_0\in\mV$. Then it holds that for any $T>0$ and $0\leq s<T$
\be
&&\lim\limits_{l\rightarrow 0}\mE\sup\limits_{v\in[s,s+l]}|\bar u^{\e,n}_{u_0}(v)-\bar u^{\e,n}_{u_0}(s)|^2_\mH=0,\label{baruenrs}\\
&&\lim\limits_{l\rightarrow 0}\mE\sup\limits_{v\in[s,s+l]}|\bar u^{0,n}_{u_0}(v)-\bar u^{0,n}_{u_0}(s)|^2_\mH=0.\label{baru0nrs}
\ee
\el
\begin{proof}
Since the proofs of (\ref{baruenrs}) and (\ref{baru0nrs}) are similar, we only prove (\ref{baruenrs}).

By the It\^o formula, it holds that for $0\leq s<v<s+l\leq T$
\ce
&&|\bar u^{\e,n}_{u_0}(v)-\bar u^{\e,n}_{u_0}(s)|^2_\mH\\
&=&2\int_s^v\<\bar u^{\e,n}_{u_0}(r)-\bar u^{\e,n}_{u_0}(s),\frac{\partial^2 \bar u^{\e,n}_{u_0}(r)}{\partial x^2}\>\dif r+2\int_s^v\<\bar u^{\e,n}_{u_0}(r)-\bar u^{\e,n}_{u_0}(s),\frac{\partial g(r,\bar u^{\e,n}_{u_0}(r))}{\partial x}\> \dif r\\
&&+2\int_s^v\<\bar u^{\e,n}_{u_0}(r)-\bar u^{\e,n}_{u_0}(s),f(\frac{r}{\e},\cdot,\bar u^{\e,n}_{u_0}(r))\>\dif r\\
&&+2\sum\limits_{j=1}^d\int_s^v\<\bar u^{\e,n}_{u_0}(r)-\bar u^{\e,n}_{u_0}(s),\sigma_j(\frac{r}{\e},\cdot,\bar u^{\e,n}_{u_0}(r))\>\dif W_j(r)\\
&&+2\int_s^v\<\bar u^{\e,n}_{u_0}(r)-\bar u^{\e,n}_{u_0}(s),n\bar u^{\e,n}_{u_0}(r)^-\>\dif r+\sum\limits_{j=1}^d\int_s^v|\sigma_j(\frac{r}{\e},\cdot,\bar u^{\e,n}_{u_0}(r))|^2_\mH\dif r\\
&=:&\cI_1+\cI_2+\cI_3+\cI_4+\cI_5+\cI_6.
\de

For $\cI_1$, the integration by parts and the Young inequality imply that
\ce
\cI_1&=&-2\int_s^v\<\frac{\p(\bar u^{\e,n}_{u_0}(r)-\bar u^{\e,n}_{u_0}(s))}{\p x},\frac{\partial \bar u^{\e,n}_{u_0}(r)}{\partial x}\>\dif r\\
&=&-2\int_s^v\<\frac{\p(\bar u^{\e,n}_{u_0}(r)-\bar u^{\e,n}_{u_0}(s))}{\p x},\frac{\p(\bar u^{\e,n}_{u_0}(r)-\bar u^{\e,n}_{u_0}(s))}{\p x}\>\dif r\\
&&-2\int_s^v\<\frac{\p(\bar u^{\e,n}_{u_0}(r)-\bar u^{\e,n}_{u_0}(s))}{\p x},\frac{\partial \bar u^{\e,n}_{u_0}(s)}{\partial x}\>\dif r\\
&\leq&-2\int_s^v\|\bar u^{\e,n}_{u_0}(r)-\bar u^{\e,n}_{u_0}(s)\|^2_\mV\dif r+2\int_s^v\left|\frac{\p(\bar u^{\e,n}_{u_0}(r)-\bar u^{\e,n}_{u_0}(s))}{\p x}\right|_\mH\left|\frac{\partial \bar u^{\e,n}_{u_0}(s)}{\partial x}\right|_\mH\dif r\\
&\leq&-2\int_s^v\|\bar u^{\e,n}_{u_0}(r)-\bar u^{\e,n}_{u_0}(s)\|^2_\mV\dif r+\int_s^v\|\bar u^{\e,n}_{u_0}(r)-\bar u^{\e,n}_{u_0}(s)\|^2_\mV\dif r+\int_s^v\|\bar u^{\e,n}_{u_0}(s)\|^2_\mV\dif r\\
&=&-\int_s^v\|\bar u^{\e,n}_{u_0}(r)-\bar u^{\e,n}_{u_0}(s)\|^2_\mV\dif r+\|\bar u^{\e,n}_{u_0}(s)\|^2_\mV(v-s).
\de
For $\cI_2$, from the integration by parts, the Young inequality and $(\mathbf{H}^C_{g})$, it follows that
\ce
\cI_2&=&-2\int_s^v\<\frac{\p(\bar u^{\e,n}_{u_0}(r)-\bar u^{\e,n}_{u_0}(s))}{\p x},g(r,\bar u^{\e,n}_{u_0}(r))\>\dif r\\
&\leq&2\int_s^v\left|\frac{\p(\bar u^{\e,n}_{u_0}(r)-\bar u^{\e,n}_{u_0}(s))}{\p x}\right|_\mH\left|g(r,\bar u^{\e,n}_{u_0}(r))\right|_\mH\dif r\\
&\leq&\int_s^v\|\bar u^{\e,n}_{u_0}(r)-\bar u^{\e,n}_{u_0}(s)\|^2_\mV\dif r\\
&&+\int_s^v\int_0^12\(|g(r,0)|^2+2L_g^2(1+|\bar u^{\e,n}_{u_0}(r,x)|^2)|\bar u^{\e,n}_{u_0}(r,x)|^2\)\dif x\dif r\\
&\leq&\int_s^v\|\bar u^{\e,n}_{u_0}(r)-\bar u^{\e,n}_{u_0}(s)\|^2_\mV\dif r+2\sup\limits_{r\in[0,T]}|g(r,0)|^2(v-s)\\
&&+\int_s^v|\bar u^{\e,n}_{u_0}(r)|^2_{L^\infty([0,1])}4L_g^2(1+|\bar u^{\e,n}_{u_0}(r)|^2_\mH)\dif r\\
&\leq&\int_s^v\|\bar u^{\e,n}_{u_0}(r)-\bar u^{\e,n}_{u_0}(s)\|^2_\mV\dif r+2\sup\limits_{r\in[0,T]}|g(r,0)|^2(v-s)\\
&&+4L_g^2(1+\sup\limits_{r\in[0,T]}|\bar u^{\e,n}_{u_0}(r)|^2_\mH)\int_s^v\|\bar u^{\e,n}_{u_0}(r)\|^2_\mV\dif r,
\de
where we use the fact $|\bar u^{\e,n}_{u_0}(r)|_{L^\infty([0,1])}\leq \|\bar u^{\e,n}_{u_0}(r)\|_\mV$. For $\cI_3$, by $(\mathbf{H}_{f})$, we infer that
\ce
\cI_3\leq \int_s^v|\bar u^{\e,n}_{u_0}(r)-\bar u^{\e,n}_{u_0}(s)|^2_\mH\dif r+L_f\int_s^v(1+|\bar u^{\e,n}_{u_0}(r)|^2_\mH)\dif r.
\de
Besides, note that $\bar u^{\e,n}_{u_0}(s,x)\geq 0$. Thus, it holds that
\ce
\cI_5&=&2\int_s^v\int_0^1(\bar u^{\e,n}_{u_0}(r,x)-\bar u^{\e,n}_{u_0}(s,x))n\bar u^{\e,n}_{u_0}(r,x)^-\dif x\dif r\\
&=&-2n\int_s^v|\bar u^{\e,n}_{u_0}(r,x)^-|^2_\mH\dif r-2n\int_s^v\int_0^1\bar u^{\e,n}_{u_0}(s,x))\bar u^{\e,n}_{u_0}(r,x)^-\dif x\dif r\\
&\leq& 0.
\de

Next, collecting the above deduction and applying the BDG inequality, we conclude that
\ce
&&\mE\sup\limits_{v\in[s,s+l]}|\bar u^{\e,n}_{u_0}(v)-\bar u^{\e,n}_{u_0}(s)|^2_\mH\\
&\leq& l\mE\|\bar u^{\e,n}_{u_0}(s)\|^2_\mV+2l\sup\limits_{r\in[0,T]}|g(r,0)|^2+4L_g^2\mE(1+\sup\limits_{r\in[0,T]}|\bar u^{\e,n}_{u_0}(r)|^2_\mH)\int_s^{s+l}\|\bar u^{\e,n}_{u_0}(r)\|^2_\mV\dif r\\
&&+\int_s^{s+l}\mE|\bar u^{\e,n}_{u_0}(r)-\bar u^{\e,n}_{u_0}(s)|^2_\mH\dif r+L_f\int_s^{s+l}(1+\mE|\bar u^{\e,n}_{u_0}(r)|^2_\mH)\dif r\\
&&+C\sum\limits_{j=1}^d\mE\left(\int_s^{s+l}|\bar u^{\e,n}_{u_0}(r)-\bar u^{\e,n}_{u_0}(s)|^2_\mH|\sigma_j(\frac{r}{\e},\cdot,\bar u^{\e,n}_{u_0}(r))|^2_\mH\dif r\right)^{1/2}\\
&&+\sum\limits_{j=1}^d\int_s^{s+l}\mE|\sigma_j(\frac{r}{\e},\cdot,\bar u^{\e,n}_{u_0}(r))|^2_\mH\dif r\\
&\leq&l\mE\|\bar u^{\e,n}_{u_0}(s)\|^2_\mV+2l\sup\limits_{r\in[0,T]}|g(r,0)|^2\\
&&+4L_g^2\left(\mE(1+\sup\limits_{r\in[0,T]}|\bar u^{\e,n}_{u_0}(r)|^2_\mH)^2\right)^{1/2}\left(\mE\left(\int_s^{s+l}\|\bar u^{\e,n}_{u_0}(r)\|^2_\mV\dif r\right)^2\right)^{1/2}\\
&&+\int_s^{s+l}\mE|\bar u^{\e,n}_{u_0}(r)-\bar u^{\e,n}_{u_0}(s)|^2_\mH\dif r+C\int_s^{s+l}(1+\mE|\bar u^{\e,n}_{u_0}(r)|^2_\mH)\dif r\\
&&+\frac{1}{2}\mE\sup\limits_{v\in[s,s+l]}|\bar u^{\e,n}_{u_0}(v)-\bar u^{\e,n}_{u_0}(s)|^2_\mH.
\de
The Gronwall inequality and (\ref{baruenes4}) imply (\ref{baruenrs}). The proof is complete.
\end{proof}

\bp\label{baruenbaru0nprob}
Assume that $(\mathbf{H}_{g})$, $(\mathbf{H}_{f})$, $(\mathbf{H}_{\s})$, $(\mathbf{H}^C_{g})$ and $(\mathbf{H}^L_{f,\s})$ hold and $u_0\in\mV$. Then it holds that for any $\d>0$
\ce
\lim\limits_{\e\rightarrow0}\mP\left\{\sup\limits_{t\in[0,T]}|\bar u^{\e,n}_{u_0}(t)-\bar u^{0,n}_{u_0}(t)|^2_\mH+\int_0^T\|\bar u^{\e,n}_{u_0}(t)-\bar u^{0,n}_{u_0}(t)\|^2_\mV\dif t\geq\d\right\}=0.
\de
\ep
\begin{proof}
{\bf Step 1.} We prove that there exists a constant $\a_4>0$ such that
\be
&&\lim\limits_{\e\rightarrow0}\Bigg\{\mE\left[\sup\limits_{t\in[0,T]}\exp\left\{-\a_4\int_0^t(1+\|\bar u^{\e,n}_{u_0}(s)\|^2_\mV+\|\bar u^{0,n}_{u_0}(s)\|^2_\mV)\dif s\right\}|\bar u^{\e,n}_{u_0}(t)-\bar u^{0,n}_{u_0}(t)|^2_\mH\right]\no\\
&&\quad\quad+\mE\left[\int_0^T\exp\left\{-\a_4\int_0^t(1+\|\bar u^{\e,n}_{u_0}(s)\|^2_\mV+\|\bar u^{0,n}_{u_0}(s)\|^2_\mV)\dif s\right\}\|\bar u^{\e,n}_{u_0}(t)-\bar u^{0,n}_{u_0}(t)\|^2_\mV\dif t\right]\Bigg\}\no\\
&&=0.
\label{baruenbaru0ncon}
\ee

First of all, for any $\a_4>0$, set 
\ce
\bar\Sigma(t):=\exp\left\{-\a_4\int_0^t(1+\|\bar u^{\e,n}_{u_0}(s)\|^2_\mV+\|\bar u^{0,n}_{u_0}(s)\|^2_\mV)\dif s\right\},
\de
and by the It\^o formula and integration by parts, it holds that
\ce
&&\bar\Sigma(t)|\bar u^{\e,n}_{u_0}(t)-\bar u^{0,n}_{u_0}(t)|^2_\mH+2\int_0^t\bar\Sigma(s)\|\bar u^{\e,n}_{u_0}(t)-\bar u^{0,n}_{u_0}(t)\|^2_\mV\dif s\no\\
&=&-\a_4\int_0^t\bar\Sigma(s)(1+\|\bar u^{\e,n}_{u_0}(s)\|^2_\mV+\|\bar u^{0,n}_{u_0}(s)\|^2_\mV)|\bar u^{\e,n}_{u_0}(s)-\bar u^{0,n}_{u_0}(s)|^2_\mH\dif s\no\\
&&+2\int_0^t\bar\Sigma(s)\<\bar u^{\e,n}_{u_0}(s)-\bar u^{0,n}_{u_0}(s), \frac{\partial \(g(s,\bar u^{\e,n}_{u_0}(s))-g(s,\bar u^{0,n}_{u_0}(s))\)}{\partial x}\>\dif s\no\\
&&+2\int_0^t\bar\Sigma(s)\<\bar u^{\e,n}_{u_0}(s)-\bar u^{0,n}_{u_0}(s), f(\frac{s}{\e},\cdot,\bar u^{\e,n}_{u_0}(s))-\bar f(\cdot,\bar u^{0,n}_{u_0}(s))\>\dif s\no\\
&&+2\sum\limits_{j=1}^d\int_0^t\bar\Sigma(s)\<\bar u^{\e,n}_{u_0}(s)-\bar u^{0,n}_{u_0}(s), \sigma_j(\frac{s}{\e},\cdot,\bar u^{\e,n}_{u_0}(s))-\bar\sigma_j(\cdot,\bar u^{0,n}_{u_0}(s))\>\dif W_j(s)\no\\
&&+2n\int_0^t\bar\Sigma(s)\<\bar u^{\e,n}_{u_0}(s)-\bar u^{0,n}_{u_0}(s),\bar u^{\e,n}_{u_0}(s)^- -\bar u^{0,n}_{u_0}(s)^-\>\dif s\no\\
&&+\sum\limits_{j=1}^d\int_0^t\bar\Sigma(s)|\sigma_j(\frac{s}{\e},\cdot,\bar u^{\e,n}_{u_0}(s))-\bar\sigma_j(\cdot,\bar u^{0,n}_{u_0}(s))|^2_\mH\dif s\no\\
&=:&\cJ_1+\cJ_2+\cJ_3+\cJ_4+\cJ_5+\cJ_6.
\de
For $\cJ_2$ and $\cJ_5$, by the similar deduction to that for (\ref{i2}) and (\ref{j6}), we obtain that
\ce
\cJ_2&\leq&\int_0^t\bar\Sigma(s)\|\bar u^{\e,n}_{u_0}(s)-\bar u^{0,n}_{u_0}(s)\|^2_\mV\dif s\no\\
&&+3L_g^2\int_0^t\bar\Sigma(s)(1+\|\bar u^{\e,n}_{u_0}(s)\|^2_\mV+\|\bar u^{0,n}_{u_0}(s)\|^2_\mV)|\bar u^{\e,n}_{u_0}(s)-\bar u^{0,n}_{u_0}(s)|^2_\mH\dif s,\\
\cJ_5&\leq&2n\int_0^t\bar\Sigma(s)|\bar u^{\e,n}_{u_0}(s)-\bar u^{0,n}_{u_0}(s)|^2_\mH\dif s.
\de
For $\cJ_3$, $(\mathbf{H}_{f})$ and $(\mathbf{H}^L_{f,\s})$ imply that
\ce
\cJ_3&=&2\int_0^t\bar\Sigma(s)\<\bar u^{\e,n}_{u_0}(s)-\bar u^{0,n}_{u_0}(s), f(\frac{s}{\e},\cdot,\bar u^{\e,n}_{u_0}(s))-\bar f(\cdot,\bar u^{\e,n}_{u_0}(s))\>\dif s\no\\
&&+2\int_0^t\bar\Sigma(s)\<\bar u^{\e,n}_{u_0}(s)-\bar u^{0,n}_{u_0}(s), \bar f(\cdot,\bar u^{\e,n}_{u_0}(s))-\bar f(\cdot,\bar u^{0,n}_{u_0}(s))\>\dif s\no\\
&\leq&2\left|\int_0^t\bar\Sigma(s)\<\bar u^{\e,n}_{u_0}(s)-\bar u^{0,n}_{u_0}(s), f(\frac{s}{\e},\cdot,\bar u^{\e,n}_{u_0}(s))-\bar f(\cdot,\bar u^{\e,n}_{u_0}(s))\>\dif s\right|\no\\
&&+2L_f\int_0^t\bar\Sigma(s)|\bar u^{\e,n}_{u_0}(s)-\bar u^{0,n}_{u_0}(s)|^2_\mH\dif s.
\de

Combining the above deduction and taking $\a_4\geq 3L_g^2$, we infer that
\ce
&&\bar\Sigma(t)|\bar u^{\e,n}_{u_0}(t)-\bar u^{0,n}_{u_0}(t)|^2_\mH+\int_0^t\bar\Sigma(s)\|\bar u^{\e,n}_{u_0}(t)-\bar u^{0,n}_{u_0}(t)\|^2_\mV\dif s\\
&\leq&(2n+2L_f)\int_0^t\bar\Sigma(s)|\bar u^{\e,n}_{u_0}(s)-\bar u^{0,n}_{u_0}(s)|^2_\mH\dif s\\
&&+2\left|\int_0^t\bar\Sigma(s)\<\bar u^{\e,n}_{u_0}(s)-\bar u^{0,n}_{u_0}(s), f(\frac{s}{\e},\cdot,\bar u^{\e,n}_{u_0}(s))-\bar f(\cdot,\bar u^{\e,n}_{u_0}(s))\>\dif s\right|\\
&&+|\cJ_4|+\cJ_{6},
\de
and furthermore
\be
&&\mE\left[\sup\limits_{t\in[0,T]}\bar\Sigma(t)|\bar u^{\e,n}_{u_0}(t)-\bar u^{0,n}_{u_0}(t)|^2_\mH+\int_0^T\bar\Sigma(s)\|\bar u^{\e,n}_{u_0}(t)-\bar u^{0,n}_{u_0}(t)\|^2_\mV\dif s\right]\no\\
&\leq&(2n+2L_f)\int_0^T\mE\bar\Sigma(s)|\bar u^{\e,n}_{u_0}(s)-\bar u^{0,n}_{u_0}(s)|^2_\mH\dif s\no\\
&&+2\mE\sup\limits_{t\in[0,T]}\left|\int_0^t\bar\Sigma(s)\<\bar u^{\e,n}_{u_0}(s)-\bar u^{0,n}_{u_0}(s), f(\frac{s}{\e},\cdot,\bar u^{\e,n}_{u_0}(s))-\bar f(\cdot,\bar u^{\e,n}_{u_0}(s))\>\dif s\right|\no\\
&&+\mE\sup\limits_{t\in[0,T]}|\cJ_4|+\mE\sup\limits_{t\in[0,T]}\cJ_{6}.
\label{cj35}
\ee
For $\mE\sup\limits_{t\in[0,T]}|\cJ_4|+\mE\sup\limits_{t\in[0,T]}\cJ_{6}$, the BDG inequality and the Young inequality yield that
\be
&&\mE\sup\limits_{t\in[0,T]}|\cJ_4|+\mE\sup\limits_{t\in[0,T]}\cJ_{6}\no\\
&\leq& C\sum\limits_{j=1}^d\mE\left(\int_0^T\bar\Sigma(s)^2|\bar u^{\e,n}_{u_0}(s)-\bar u^{0,n}_{u_0}(s)|^2_\mH|\sigma_j(\frac{s}{\e},\cdot,\bar u^{\e,n}_{u_0}(s))-\bar\sigma_j(\cdot,\bar u^{0,n}_{u_0}(s))|^2_\mH\dif s\right)^{1/2}\no\\
&&+\sum\limits_{j=1}^d\mE\int_0^T\bar\Sigma(s)|\sigma_j(\frac{s}{\e},\cdot,\bar u^{\e,n}_{u_0}(s))-\bar\sigma_j(\cdot,\bar u^{0,n}_{u_0}(s))|^2_\mH\dif s\no\\
&\leq&\frac{1}{2}\mE\sup\limits_{t\in[0,T]}\bar\Sigma(t)|\bar u^{\e,n}_{u_0}(t)-\bar u^{0,n}_{u_0}(t)|^2_\mH+C\sum\limits_{j=1}^d\mE\int_0^T\bar\Sigma(s)|\sigma_j(\frac{s}{\e},\cdot,\bar u^{\e,n}_{u_0}(s))-\bar\sigma_j(\cdot,\bar u^{0,n}_{u_0}(s))|^2_\mH\dif s\no\\
&\leq&\frac{1}{2}\mE\sup\limits_{t\in[0,T]}\bar\Sigma(t)|\bar u^{\e,n}_{u_0}(t)-\bar u^{0,n}_{u_0}(t)|^2_\mH+2C\sum\limits_{j=1}^d\mE\int_0^T\bar\Sigma(s)|\sigma_j(\frac{s}{\e},\cdot,\bar u^{\e,n}_{u_0}(s))-\bar\sigma_j(\cdot,\bar u^{\e,n}_{u_0}(s))|^2_\mH\dif s\no\\
&&+2C\sum\limits_{j=1}^d\mE\int_0^t\bar\Sigma(s)|\bar\sigma_j(\cdot,\bar u^{\e,n}_{u_0}(s))-\bar\sigma_j(\cdot,\bar u^{0,n}_{u_0}(s))|^2_\mH\dif s\no\\
&\leq&\frac{1}{2}\mE\sup\limits_{t\in[0,T]}\bar\Sigma(t)|\bar u^{\e,n}_{u_0}(t)-\bar u^{0,n}_{u_0}(t)|^2_\mH+2C\sum\limits_{j=1}^d\mE\int_0^T\bar\Sigma(s)|\sigma_j(\frac{s}{\e},\cdot,\bar u^{\e,n}_{u_0}(s))-\bar\sigma_j(\cdot,\bar u^{\e,n}_{u_0}(s))|^2_\mH\dif s\no\\
&&+CL_\s\int_0^T\mE\bar\Sigma(s)|\bar u^{\e,n}_{u_0}(s)-\bar u^{0,n}_{u_0}(s)|^2_\mH\dif s.
\label{cj5}
\ee
Collecting (\ref{cj35}) and (\ref{cj5}) and applying the Gronwall inequality, we conclude that
\be
&&\mE\left[\sup\limits_{t\in[0,T]}\bar\Sigma(t)|\bar u^{\e,n}_{u_0}(t)-\bar u^{0,n}_{u_0}(t)|^2_\mH+\int_0^T\bar\Sigma(s)\|\bar u^{\e,n}_{u_0}(t)-\bar u^{0,n}_{u_0}(t)\|^2_\mV\dif s\right]\no\\
&\leq&C\mE\sup\limits_{t\in[0,T]}\left|\int_0^t\bar\Sigma(s)\<\bar u^{\e,n}_{u_0}(s)-\bar u^{0,n}_{u_0}(s), f(\frac{s}{\e},\cdot,\bar u^{\e,n}_{u_0}(s))-\bar f(\cdot,\bar u^{\e,n}_{u_0}(s))\>\dif s\right|\no\\
&&+2C\sum\limits_{j=1}^d\mE\int_0^T\bar\Sigma(s)|\sigma_j(\frac{s}{\e},\cdot,\bar u^{\e,n}_{u_0}(s))-\bar\sigma_j(\cdot,\bar u^{\e,n}_{u_0}(s))|^2_\mH\dif s\no\\
&=:&\cZ_1+\cZ_2.
\label{grocj3}
\ee 

Finally, note that
\be
&&\lim\limits_{\e\rightarrow 0}\cZ_1=0,\label{cj3}\\
&&\lim\limits_{\e\rightarrow 0}\cZ_2=0,\label{cj30}
\ee
which together with (\ref{grocj3}) yields that (\ref{baruenbaru0ncon}) holds. Thus, we complete the proof of (\ref{baruenbaru0ncon}). 

{\bf Step 2.} We prove (\ref{cj3}) and (\ref{cj30}).

Note that for any $0<\t<T$,
\ce
&&\<\bar u^{\e,n}_{u_0}(s)-\bar u^{0,n}_{u_0}(s), f(\frac{s}{\e},\cdot,\bar u^{\e,n}_{u_0}(s))-\bar f(\cdot,\bar u^{\e,n}_{u_0}(s))\>\\
&=&\<\bar u^{\e,n}_{u_0}(s)-\bar u^{0,n}_{u_0}(s)-\bar u^{\e,n}_{u_0}(s(\t))+\bar u^{0,n}_{u_0}(s(\t)), f(\frac{s}{\e},\cdot,\bar u^{\e,n}_{u_0}(s))-\bar f(\cdot,\bar u^{\e,n}_{u_0}(s))\>\\
&&+\<\bar u^{\e,n}_{u_0}(s(\t))-\bar u^{0,n}_{u_0}(s(\t)), f(\frac{s}{\e},\cdot,\bar u^{\e,n}_{u_0}(s))-f(\frac{s}{\e},\cdot,\bar u^{\e,n}_{u_0}(s(\t)))\>\\
&&+\<\bar u^{\e,n}_{u_0}(s(\t))-\bar u^{0,n}_{u_0}(s(\t)), \bar f(\cdot,\bar u^{\e,n}_{u_0}(s(\t)))-\bar f(\cdot,\bar u^{\e,n}_{u_0}(s))\>\\
&&+\<\bar u^{\e,n}_{u_0}(s(\t))-\bar u^{0,n}_{u_0}(s(\t)), f(\frac{s}{\e},\cdot,\bar u^{\e,n}_{u_0}(s(\t)))-\bar f(\cdot,\bar u^{\e,n}_{u_0}(s(\t)))\>,
\de
where $s(\t):=[\frac{s}{\t}]\t$ and $[\frac{s}{\t}]$ denotes the integer part of $\frac{s}{\t}$. Thus, it holds that
\ce
\cZ_1&\leq&\mE\int_0^T|\bar u^{\e,n}_{u_0}(s)-\bar u^{0,n}_{u_0}(s)-\bar u^{\e,n}_{u_0}(s(\t))+\bar u^{0,n}_{u_0}(s(\t))|_\mH|f(\frac{s}{\e},\cdot,\bar u^{\e,n}_{u_0}(s))-\bar f(\cdot,\bar u^{\e,n}_{u_0}(s))|_\mH\dif s\\
&&+\mE\int_0^T|\bar u^{\e,n}_{u_0}(s(\t))-\bar u^{0,n}_{u_0}(s(\t))|_\mH|f(\frac{s}{\e},\cdot,\bar u^{\e,n}_{u_0}(s))-f(\frac{s}{\e},\cdot,\bar u^{\e,n}_{u_0}(s(\t)))|_\mH\dif s\\
&&+\mE\int_0^T|\bar u^{\e,n}_{u_0}(s(\t))-\bar u^{0,n}_{u_0}(s(\t))|_\mH|\bar f(\cdot,\bar u^{\e,n}_{u_0}(s(\t)))-\bar f(\cdot,\bar u^{\e,n}_{u_0}(s))|_\mH\dif s\\
&&+\mE\int_0^T|\bar u^{\e,n}_{u_0}(s(\t))-\bar u^{0,n}_{u_0}(s(\t))|_\mH|f(\frac{s}{\e},\cdot,\bar u^{\e,n}_{u_0}(s(\t)))-\bar f(\cdot,\bar u^{\e,n}_{u_0}(s(\t)))|_\mH\dif s\\
&=:&\cZ_{11}+\cZ_{12}+\cZ_{13}+\cZ_{14},
\de
where we use the fact that $\bar\Sigma(s)\leq 1$.

For $\cZ_{11}$, the H\"older inequality, $(\mathbf{H}_{f})$ and $(\mathbf{H}^L_{f,\s})$ imply that
\ce
\cZ_{11}&\leq&\left(\mE\int_0^T|\bar u^{\e,n}_{u_0}(s)-\bar u^{0,n}_{u_0}(s)-\bar u^{\e,n}_{u_0}(s(\t))+\bar u^{0,n}_{u_0}(s(\t))|^2_\mH\dif s\right)^{1/2}\\
&&\times\left(\mE\int_0^T|f(\frac{s}{\e},\cdot,\bar u^{\e,n}_{u_0}(s))-\bar f(\cdot,\bar u^{\e,n}_{u_0}(s))|^2_\mH\dif s\right)^{1/2}\\
&\leq&\left(2\mE\int_0^T\(|\bar u^{\e,n}_{u_0}(s)-\bar u^{\e,n}_{u_0}(s(\t))|^2_\mH+|\bar u^{0,n}_{u_0}(s)-\bar u^{0,n}_{u_0}(s(\t))|^2_\mH\)\dif s\right)^{1/2}\\
&&\times\left(2\mE\int_0^T\(|f(\frac{s}{\e},\cdot,\bar u^{\e,n}_{u_0}(s))|^2_\mH+|\bar f(\cdot,\bar u^{\e,n}_{u_0}(s))|^2_\mH\)\dif s\right)^{1/2}\\
&\leq&\left(2\int_0^T\(\mE\sup\limits_{v\in[s,s+\t]}|\bar u^{\e,n}_{u_0}(v)-\bar u^{\e,n}_{u_0}(s)|^2_\mH+\mE\sup\limits_{v\in[s,s+\t]}|\bar u^{0,n}_{u_0}(v)-\bar u^{0,n}_{u_0}(s)|^2_\mH\)\dif s\right)^{1/2}\\
&&\times C(1+|u_0|_\mH).
\de
For $\cZ_{12}$ and $\cZ_{13}$, by the H\"older inequality, $(\mathbf{H}_{f})$ and $(\mathbf{H}^L_{f,\s})$, we infer that
\ce
\cZ_{12}+\cZ_{13}&\leq& C\left(\mE\int_0^T|\bar u^{\e,n}_{u_0}(s(\t))-\bar u^{0,n}_{u_0}(s(\t))|^2_\mH\dif s\right)^{1/2}\\
&&\times\left(\mE\int_0^T|\bar u^{\e,n}_{u_0}(s)-\bar u^{\e,n}_{u_0}(s(\t))|^2_\mH\dif s\right)^{1/2}\\
&\leq& C\left(\int_0^T2\(\mE\sup\limits_{s\in[0,T]}|\bar u^{\e,n}_{u_0}(s)|^2_\mH+\mE\sup\limits_{s\in[0,T]}|\bar u^{0,n}_{u_0}(s)|^2_\mH\)\dif s\right)^{1/2}\\
&&\times\left(\int_0^T\mE\sup\limits_{v\in[s,s+\t]}|\bar u^{\e,n}_{u_0}(v)-\bar u^{\e,n}_{u_0}(s)|^2_\mH\dif s\right)^{1/2}\\
&\leq&C(1+|u_0|_\mH)\left(\int_0^T\mE\sup\limits_{v\in[s,s+\t]}|\bar u^{\e,n}_{u_0}(v)-\bar u^{\e,n}_{u_0}(s)|^2_\mH\dif s\right)^{1/2}.
\de

For $\cZ_{14}$, we divide it into two parts. That is,
\ce
\cZ_{14}&=&\mE\int_0^{[T/\t]\t}|\bar u^{\e,n}_{u_0}(s(\t))-\bar u^{0,n}_{u_0}(s(\t))|_\mH|f(\frac{s}{\e},\cdot,\bar u^{\e,n}_{u_0}(s(\t)))-\bar f(\cdot,\bar u^{\e,n}_{u_0}(s(\t)))|_\mH\dif s\\
&&+\mE\int_{[T/\t]\t}^T|\bar u^{\e,n}_{u_0}(s(\t))-\bar u^{0,n}_{u_0}(s(\t))|_\mH|f(\frac{s}{\e},\cdot,\bar u^{\e,n}_{u_0}(s(\t)))-\bar f(\cdot,\bar u^{\e,n}_{u_0}(s(\t)))|_\mH\dif s\\
&=:&\cZ_{141}+\cZ_{142}.
\de
On the one hand, for $\cZ_{141}$, $(\mathbf{H}_{f})$ and $(\mathbf{H}^L_{f,\s})$ imply that
\ce
\cZ_{141}&=&\mE\sum_{k=0}^{[T/\t]-1}\int_{k\t}^{(k+1)\t}|\bar u^{\e,n}_{u_0}(s(\t))-\bar u^{0,n}_{u_0}(s(\t))|_\mH|f(\frac{s}{\e},\cdot,\bar u^{\e,n}_{u_0}(s(\t)))-\bar f(\cdot,\bar u^{\e,n}_{u_0}(s(\t)))|_\mH\dif s\\
&=&\mE\sum_{k=0}^{[T/\t]-1}\int_{k\t}^{(k+1)\t}|\bar u^{\e,n}_{u_0}(k\t)-\bar u^{0,n}_{u_0}(k\t)|_\mH|f(\frac{s}{\e},\cdot,\bar u^{\e,n}_{u_0}(k\t))-\bar f(\cdot,\bar u^{\e,n}_{u_0}(k\t))|_\mH\dif s\\
&\leq&\t\mE\sum_{k=0}^{[T/\t]-1}|\bar u^{\e,n}_{u_0}(k\t)-\bar u^{0,n}_{u_0}(k\t)|_\mH \frac{\e}{\t}\int_{\frac{k\t}{\e}}^{\frac{(k+1)\t}{\e}}|f(v,\cdot,\bar u^{\e,n}_{u_0}(k\t))-\bar f(\cdot,\bar u^{\e,n}_{u_0}(k\t))|_\mH\dif v\\
&\leq&\t\sum_{k=0}^{[T/\t]-1}\left(\mE|\bar u^{\e,n}_{u_0}(k\t)-\bar u^{0,n}_{u_0}(k\t)|^2_\mH\right)^{1/2}\\
&&\times\left(\mE\frac{\e}{\t}\int_{\frac{k\t}{\e}}^{\frac{(k+1)\t}{\e}}|f(v,\cdot,\bar u^{\e,n}_{u_0}(k\t))-\bar f(\cdot,\bar u^{\e,n}_{u_0}(k\t))|^2_\mH\dif v\right)^{1/2}\\
&\leq&\t\sum_{k=0}^{[T/\t]-1}\left(2\mE\sup\limits_{s\in[0,T]}|\bar u^{\e,n}_{u_0}(s)|^2_\mH+2\mE\sup\limits_{s\in[0,T]}|\bar u^{0,n}_{u_0}(s)|^2_\mH\right)^{1/2}\left(\mE\kappa(\frac{\t}{\e})(1+|\bar u^{\e,n}_{u_0}(k\t)|_\mH^2)\right)^{1/2}\\
&\leq&CT(1+|u_0|^2_\mH)\kappa^{1/2}(\frac{\t}{\e}).
\de
On the other hand, for $\cZ_{142}$, by $(\mathbf{H}_{f})$ and $(\mathbf{H}^L_{f,\s})$, it holds that
\ce
\cZ_{142}&\leq&\left(\mE\int_{[T/\t]\t}^T|\bar u^{\e,n}_{u_0}(s(\t))-\bar u^{0,n}_{u_0}(s(\t))|^2_\mH\dif s\right)^{1/2}\\
&&\times\left(\mE\int_{[T/\t]\t}^T|f(\frac{s}{\e},\cdot,\bar u^{\e,n}_{u_0}(s(\t)))-\bar f(\cdot,\bar u^{\e,n}_{u_0}(s(\t)))|^2_\mH\dif s\right)^{1/2}\\
&\leq&\left(\int_{[T/\t]\t}^T2(\mE\sup\limits_{s\in[0,T]}|\bar u^{\e,n}_{u_0}(s)|^2_\mH+\mE\sup\limits_{s\in[0,T]}|\bar u^{0,n}_{u_0}(s)|^2_\mH)\dif s\right)^{1/2}\\
&&\times\left(C\int_{[T/\t]\t}^T(1+\mE\sup\limits_{s\in[0,T]}|\bar u^{\e,n}_{u_0}(s)|^2_\mH)\dif s\right)^{1/2}\\
&\leq&C\t(1+|u_0|^2_\mH).
\de

At last, collecting the above deduction, we conclude that
\ce
\cZ_1&\leq&C(1+|u_0|_\mH)\Bigg(\int_0^T\(\mE\sup\limits_{v\in[s,s+\t]}|\bar u^{\e,n}_{u_0}(v)-\bar u^{\e,n}_{u_0}(s)|^2_\mH\\
&&\qquad +\mE\sup\limits_{v\in[s,s+\t]}|\bar u^{0,n}_{u_0}(v)-\bar u^{0,n}_{u_0}(s)|^2_\mH\)\dif s\Bigg)^{1/2}\\
&&+C(1+|u_0|_\mH)\left(\int_0^T\mE\sup\limits_{v\in[s,s+\t]}|\bar u^{\e,n}_{u_0}(v)-\bar u^{\e,n}_{u_0}(s)|^2_\mH\dif s\right)^{1/2}\\
&&+CT(1+|u_0|^2_\mH)\kappa^{1/2}(\frac{\t}{\e})+C\t(1+|u_0|^2_\mH).
\de
Letting $\e\rightarrow 0$ first and then $\t\rightarrow 0$, by (\ref{baruenrs}) and (\ref{baru0nrs}), we have (\ref{cj3}). By the similar deduction to that for (\ref{cj3}), one can show (\ref{cj30}).

{\bf Step 3.} We prove that
\ce
\lim\limits_{\e\rightarrow0}\mP\left\{\sup\limits_{t\in[0,T]}|\bar u^{\e,n}_{u_0}(t)-\bar u^{0,n}_{u_0}(t)|^2_\mH+\int_0^T\|\bar u^{\e,n}_{u_0}(t)-\bar u^{0,n}_{u_0}(t)\|^2_\mV\dif t\geq\d\right\}=0.
\de

For any $M>0$, by the Chebyshev inequality, it holds that
\ce
&&\mP\left\{\sup\limits_{t\in[0,T]}|\bar u^{\e,n}_{u_0}(t)-\bar u^{0,n}_{u_0}(t)|^2_\mH+\int_0^T\|\bar u^{\e,n}_{u_0}(t)-\bar u^{0,n}_{u_0}(t)\|^2_\mV\dif t\geq\d\right\}\\
&=&\mP\Bigg\{\sup\limits_{t\in[0,T]}|\bar u^{\e,n}_{u_0}(t)-\bar u^{0,n}_{u_0}(t)|^2_\mH+\int_0^T\|\bar u^{\e,n}_{u_0}(t)-\bar u^{0,n}_{u_0}(t)\|^2_\mV\dif t\geq\d,\\
&&\qquad \a_4\int_0^T(1+\|\bar u^{\e,n}_{u_0}(s)\|^2_\mV+\|\bar u^{0,n}_{u_0}(s)\|^2_\mV)\dif s\leq M\Bigg\}\\
&&+\mP\Bigg\{\sup\limits_{t\in[0,T]}|\bar u^{\e,n}_{u_0}(t)-\bar u^{0,n}_{u_0}(t)|^2_\mH+\int_0^T\|\bar u^{\e,n}_{u_0}(t)-\bar u^{0,n}_{u_0}(t)\|^2_\mV\dif t\geq\d,\\
&&\qquad\a_4\int_0^T(1+\|\bar u^{\e,n}_{u_0}(s)\|^2_\mV+\|\bar u^{0,n}_{u_0}(s)\|^2_\mV)\dif s>M\Bigg\}\\
&\leq&\mP\Bigg\{\sup\limits_{t\in[0,T]}\bar\Sigma(t)|\bar u^{\e,n}_{u_0}(t)-\bar u^{0,n}_{u_0}(t)|^2_\mH+\int_0^T\bar\Sigma(t)\|\bar u^{\e,n}_{u_0}(t)-\bar u^{0,n}_{u_0}(t)\|^2_\mV\dif t\geq e^{-M}\d\Bigg\}\\
&&+\mP\left\{\a_4\int_0^T(1+\|\bar u^{\e,n}_{u_0}(s)\|^2_\mV+\|\bar u^{0,n}_{u_0}(s)\|^2_\mV)\dif s>M\right\}\\
&\leq&\frac{e^M}{\d}\mE\[\sup\limits_{t\in[0,T]}\bar\Sigma(t)|\bar u^{\e,n}_{u_0}(t)-\bar u^{0,n}_{u_0}(t)|^2_\mH+\int_0^T\bar\Sigma(t)\|\bar u^{\e,n}_{u_0}(t)-\bar u^{0,n}_{u_0}(t)\|^2_\mV\dif t\]\\
&&+\frac{\a_4}{M}\mE\int_0^T(1+\|\bar u^{\e,n}_{u_0}(s)\|^2_\mV+\|\bar u^{0,n}_{u_0}(s)\|^2_\mV)\dif s\\
&\leq&\frac{e^M}{\d}\mE\[\sup\limits_{t\in[0,T]}\bar\Sigma(t)|\bar u^{\e,n}_{u_0}(t)-\bar u^{0,n}_{u_0}(t)|^2_\mH+\int_0^T\bar\Sigma(t)\|\bar u^{\e,n}_{u_0}(t)-\bar u^{0,n}_{u_0}(t)\|^2_\mV\dif t\]\\
&&+\frac{\a_4}{M}\[T+C(1+|u_0|^2_\mH)\].
\de

Finally, letting $\e\rightarrow 0$ first and then $M\rightarrow\infty$, by (\ref{baruenbaru0ncon}) we conclude that
\ce
\lim\limits_{\e\rightarrow 0}\mP\left\{\sup\limits_{t\in[0,T]}|\bar u^{\e,n}_{u_0}(t)-\bar u^{0,n}_{u_0}(t)|^2_\mH+\int_0^T\|\bar u^{\e,n}_{u_0}(t)-\bar u^{0,n}_{u_0}(t)\|^2_\mV\dif t\geq\d\right\}=0.
\de
The proof is complete.
\end{proof}

Now, it is the position to prove Theorem \ref{averprintheo}. 

{\bf Proof of Theorem \ref{averprintheo}.}

First of all, we prove that for any $T>0$, $\bar u^{\e}_{u_0}$ converges to $\bar u^{0}_{u_0}$ in probability in $C([0,T],\mH)$. For any $\d>0$, it holds that
\ce
&&\mP\left\{\sup\limits_{t\in[0,T]}|\bar u^{\e}_{u_0}(t)-\bar u^{0}_{u_0}(t)|_\mH>\d\right\}\\
&\leq&\mP\left\{\sup\limits_{t\in[0,T]}|\bar u^{\e}_{u_0}(t)-\bar u^{\e,n}_{u_0}(t)|_\mH>\d/3\right\}+\mP\left\{\sup\limits_{t\in[0,T]}|\bar u^{\e,n}_{u_0}(t)-\bar u^{0,n}_{u_0}(t)|_\mH>\d/3\right\}\\
&&+\mP\left\{\sup\limits_{t\in[0,T]}|\bar u^{0,n}_{u_0}(t)-\bar u^{0}_{u_0}(t)|_\mH>\d/3\right\}.
\de

Next, by $(iii)$ and $(iv)$ in Lemma \ref{baruenbaru0nresu}, we know that for any $\eta>0$, there exists a $n_0\in\mN$ such that 
\ce
&&\mP\left\{\sup\limits_{t\in[0,T]}|\bar u^{\e}_{u_0}(t)-\bar u^{\e,n_0}_{u_0}(t)|_\mH>\d/3\right\}<\eta/3,\\
&&\mP\left\{\sup\limits_{t\in[0,T]}|\bar u^{0,n_0}_{u_0}(t)-\bar u^{0}_{u_0}(t)|_\mH>\d/3\right\}<\eta/3.
\de
Besides, Proposition \ref{baruenbaru0nprob} implies that there exists a $0<\e_0<1$ such that for any $\e<\e_0$
\ce
\mP\left\{\sup\limits_{t\in[0,T]}|\bar u^{\e,n_0}_{u_0}(t)-\bar u^{0,n_0}_{u_0}(t)|_\mH>\d/3\right\}<\eta/3.
\de

At last, combining the above deduction, we conclude that for any $T>0$, $\bar u^{\e}_{u_0}$ converges to $\bar u^{0}_{u_0}$ in probability in $C([0,T],\mH)$. By the same deduction to the above one, it holds that $\bar u^{\e}_{u_0}$ converges to $\bar u^{0}_{u_0}$ in probability in $L^2([0,T],\mV)$. Thus, the proof is complete.

\end{document}